\documentclass[11pt,twoside]{article}
\usepackage{amsfonts,amsmath,amssymb,amsthm}
\usepackage{epsfig}
\everymath{\displaystyle}
\def\l{\left}
\def\r{\right}

\newcommand{\bb}[1]{\begin{equation}\label{#1}}
\newcommand{\ee}{\end{equation}}
\newcommand{\bbb}{\begin{eqnarray}}
\newcommand{\eee}{\end{eqnarray}}
\newcommand{\bbbb}{\begin{eqnarray*}}
\newcommand{\eeee}{\end{eqnarray*}}
\newcommand{\nnn}{\nonumber}
\newcommand{\no}{\noindent}
\def\TT{\mathbb{T}}
\def\RR{\mathbb{R}}

\def\R#1{$(\ref{#1})$}

\begin{document}
\begin{center}

{\Large\bf An Exploration of the Approximation of Derivative
Functions via Finite Differences}

\vspace{0.2in}

Brian Jain\\
Baylor University\\
2415 South University Parks Dr., Apt 4109\\
Waco, TX 76706\\
{\sf Brian\underline{~}Jain@baylor.edu}

\vspace{0.1in}

Andrew D. Sheng\\
Westwood High School\\
12400 Mellow Meadow Drive\\
Austin, TX 78750\\
{\sf science4sail@gmail.com}

\vspace{0.2in}

\parbox[t]{5.1in}{\small{\bf Abstract: } Finite differences have
been widely used in mathematical theory as well as in scientific and
engineering computations. These concepts are constantly mentioned in
calculus. Most frequently-used difference formulas provide excellent
approximations to various derivative functions, including those used
in modeling important physical processes on uniform grids. However,
our research reveals that difference approximations on uniform grids
cannot be applied blindly on nonuniform grids, nor can difference
formulas to form consistent approximations to second derivatives. At
best, they may lose accuracy; at worst they are inconsistent.
Detailed consistency and error analysis, together with simulated
examples, will be given.}
\end{center}

\vspace{0.1in}

{\large\bf 1.
Background.~}\setcounter{section}{1}\setcounter{equation}{0} The
fundamentals of calculus were developed over a very long period of
time. They can be traced back to the late 1660s and 1670s, when
Newton and Leibniz introduced the concept of derivatives through
finite differences. The concept quickly became a cornerstone of
calculus \cite{wik,stewart1,vesely}.


\begin{center}
\epsfig{file=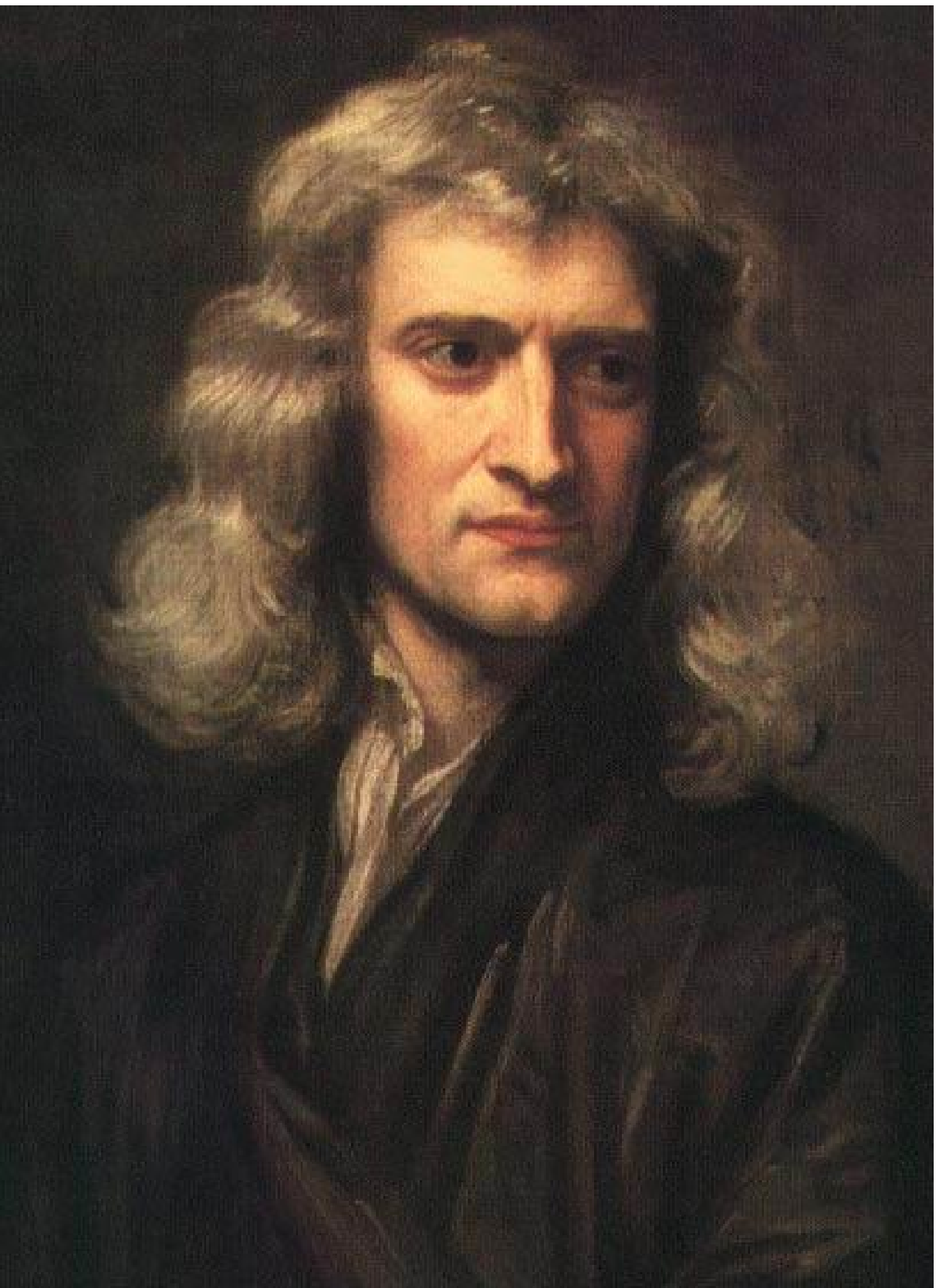,height=1.6in}~~~
\epsfig{file=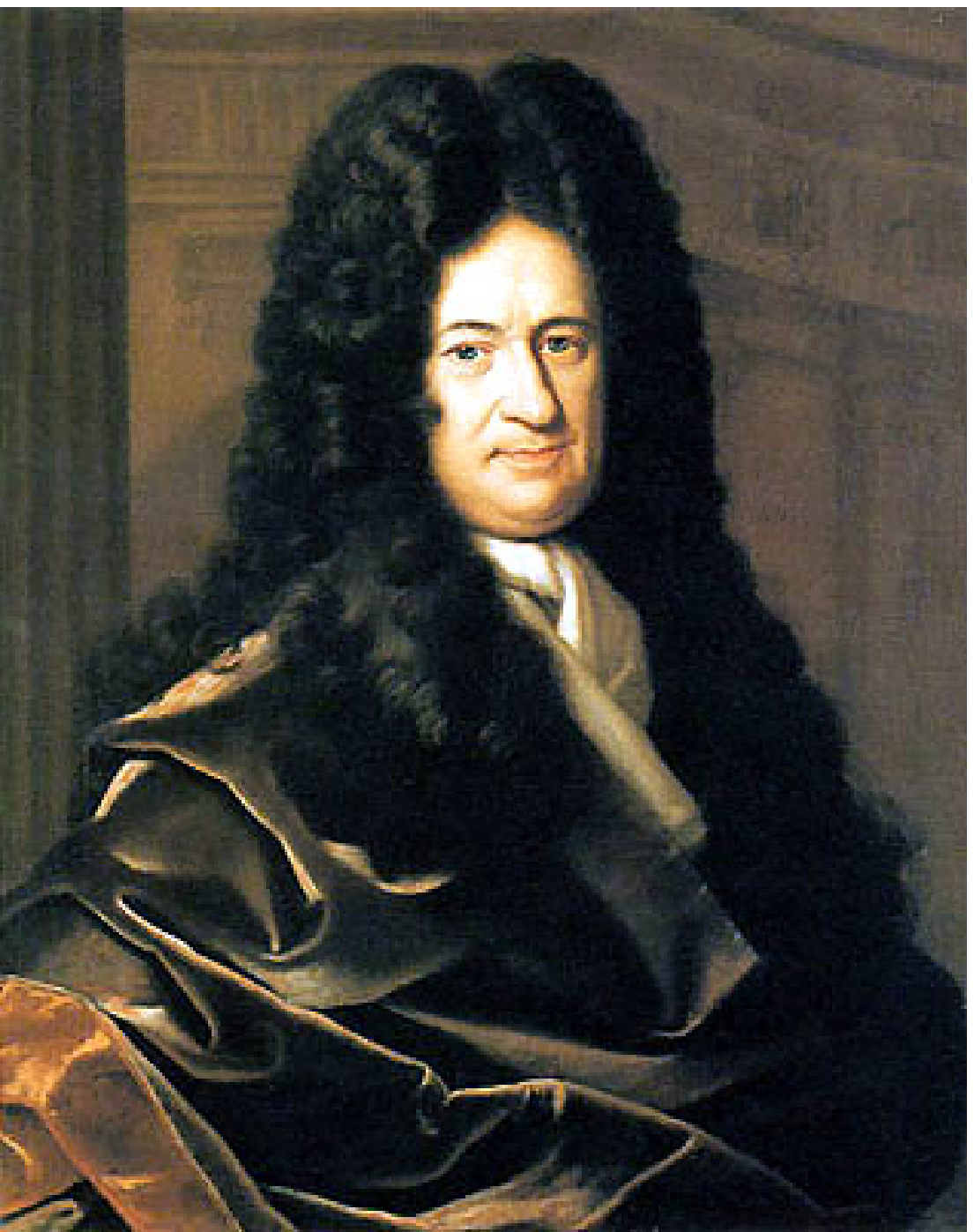,height=1.6in}

\vspace{1mm}

\parbox[t]{13.5cm}{\footnotesize{\bf Figure 1.1.} LEFT: Sir Isaac Newton
\cite{newton,vesely}. RIGHT: Gottfried Wilhelm von Leibniz
\cite{wik,leibniz}. The two scientists founded calculus as well as
the finite difference calculations.}
\end{center}

Different finite difference formulas have been constructed and
utilized for approximating the rates of change of a function in
applications \cite{wik,kelley1,levy1}. The function is given either
as an analytical expression or as a set of numbers at discrete
points in the region of interest. The mathematical models of most
science and engineering problems require such an approximation if
computer simulations are desirable \cite{iserles,rao1}.

In many discussions of one-dimensional finite differences, the
aforementioned discrete points are uniformly distributed, that is,
the distance between any two neighboring points is a positive
constant $h$ \cite{iserles,kelley1,levy1}.

Let $f(t)$ be a differentiable function on $(a, b).$ Recall the
definition of the derivative at a fixed point $t\in(a,b)$
\cite{stewart1}:
$$f'(t)=\lim_{h\rightarrow
0}\frac{f(t+h)-f(t)}{h}.$$ Thus, it is reasonable that the fraction
\bb{aa2}g(t,h)=\frac{f(t+h)-f(t)}{h},~~~0<h\ll 1,\ee would provide
an approximation of the derivative $f'(t)$ on the set $\TT=\{t,~t\pm
h,~t\pm 2h,\ldots\}\cap(a,b).$ In fact, it can be shown by
Newton-Gregory interpolation formulas that $g(t,h)$ is indeed an
approximation of $f'(t)$ \cite{vesely}. Formula \R{aa2} is called a
{\em first order forward difference\/} \cite{iserles,stewart1}.
Based on it, a {\em second order forward difference\/} could be
constructed: \bb{aa3}\frac{g(t+h,h)-g(t,h)}{h}=
\frac{f(t+2h)-2f(t+h)+f(t)}{h^2},~~~0<h\ll 1.\ee We will leave the
discussion of \R{aa3} to Section IV.

This paper studies finite difference approximations on sets of
discrete points, where distances between neighboring points vary.
Many questions remain in the situation, such as:

\begin{enumerate}
\item How good can the formula \R{aa2} be in applications? \item Are there different, or
better, approximation formulas? \item Can the formulas be used
repeatedly for approximating higher derivatives? \item How can we
evaluate the errors in approximations?
\item Can we demonstrate our study through the use of computer simulations?
\end{enumerate}

The preceding questions motivated our research in the subject. In
this article, we will study some of the most popular finite
differences for approximating the first and second order
derivatives. Cases involving different types of grids over the
intervals will be investigated. Section II will be for finite
difference preliminaries. In Section III, we will concentrate on the
first order differences. We will show that the same difference
formula may possess different accuracies on different meshes. Then,
we will show in Section IV that, although second order differences
are in general acceptable on uniform meshes, naive applications of
the formulas on nonuniform meshes will lead to incorrect results.
Computer simulated numerical examples will be given to illustrate
our conclusions in Section V. Our investigation will be summarized
in Section VI.

\vspace{0.25in}

{\large\bf 2.
Preliminaries.~}\setcounter{section}{2}\setcounter{equation}{0} Let
us consider an interval $[a, b],$ where $\infty>b>a>-\infty,$
defined in $\RR.$ We will start with definitions of the uniform and
nonuniform meshes, which are discrete sets of numbers used as
domains for our finite differences. Meshes are also called grids in
science and engineering computations \cite{rao1}.

\vspace{3mm}

\no{\bf Definition 2.1.} A {\em mesh\/} over the interval  $[a,
b]$ is defined as a set of $m+2$ distinct numbers ($m\geq 0$)
denoted as $\TT=\{t_0,~t_1,\ldots, t_{m+1}\},$ in which
$t_{k+1}>t_k,~k=0,1,\ldots,m,$ and $t_0=a,~t_{m+1}=b.$ The value
$t_k$ is called the {\em$k$th mesh point,\/} or {\em $k$th
grid,\/} of $\TT.$ Moreover, we call the quantity
$$h_k=t_{k+1}-t_k,~~~m\geq k\geq 0,$$ the {\em $k$th step size\/}
 of the
mesh. Step sizes are often less than one. $\TT$ is called a {\em
uniform mesh\/} if $h_k=h>0,~k=0,1,\ldots,m,$ otherwise $\TT$ is a
{\em nonuniform mesh.\/}

\vspace{1mm}

\begin{center}
\epsfig{file=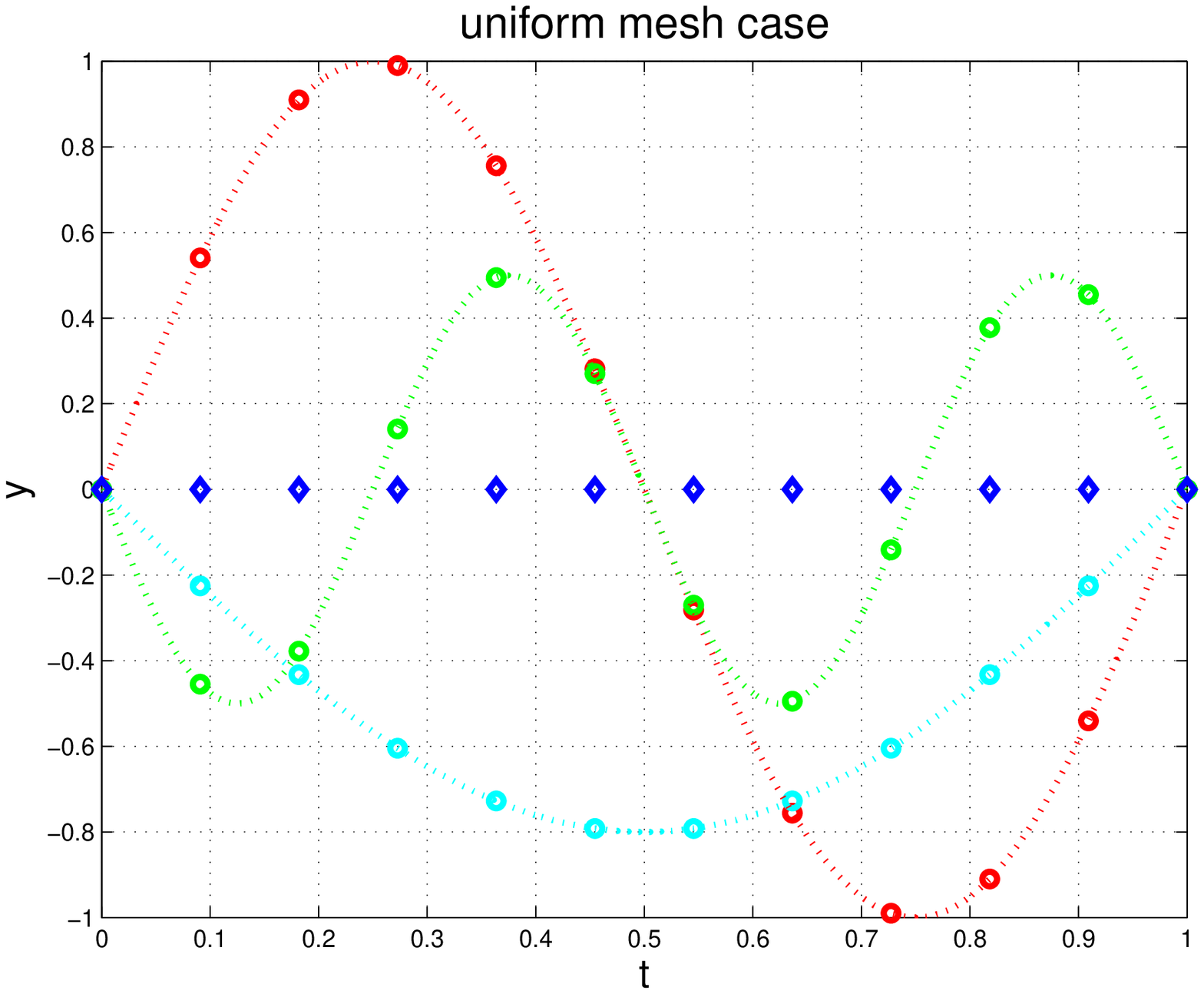,width=2.6in,height=1.6in}
\epsfig{file=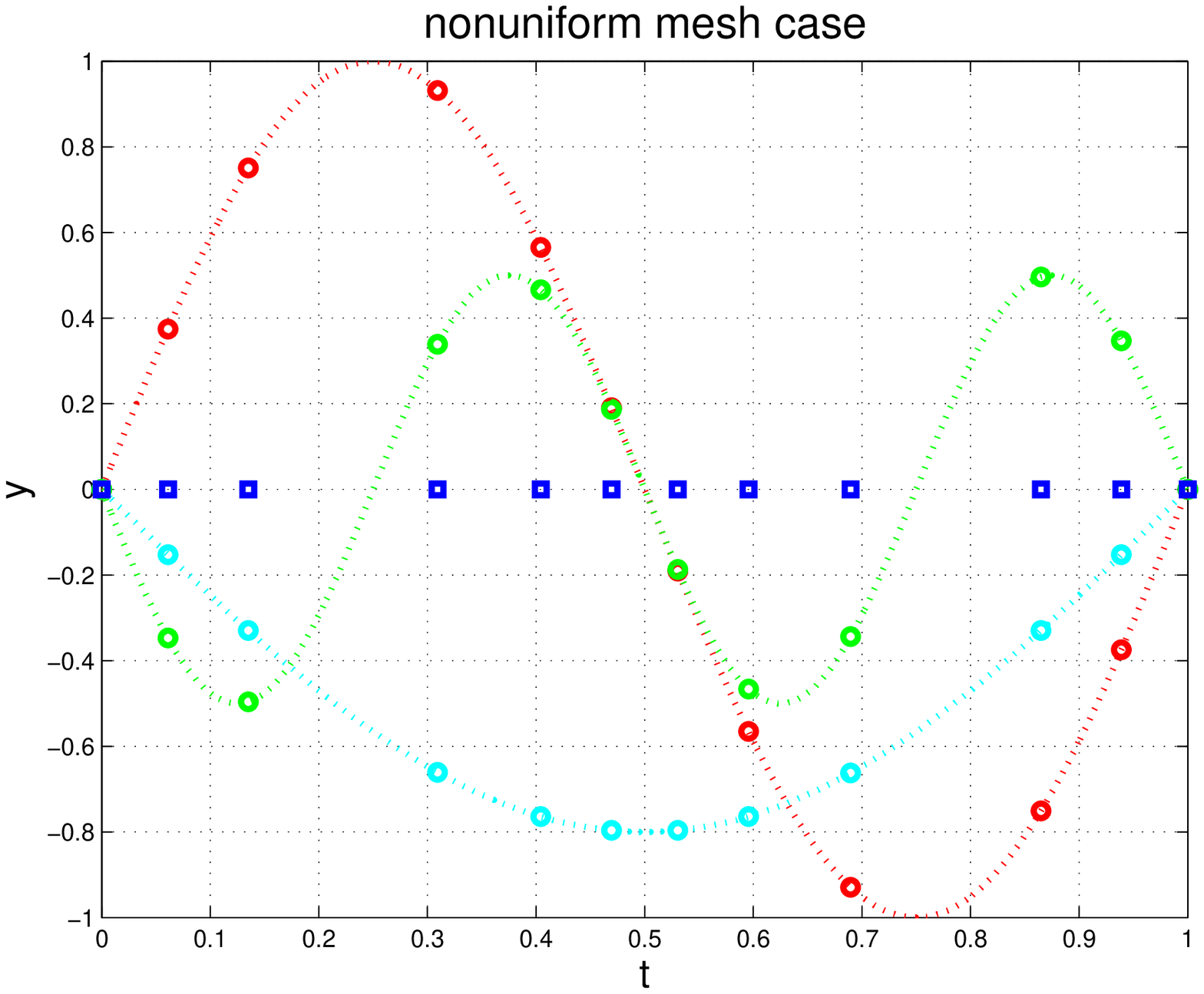,width=2.6in,height=1.6in}

\vspace{2mm}

\parbox[t]{13.5cm}{\footnotesize{\bf Figure 2.1.} Plot of the
trigonometric functions with different frequencies on the uniform
mesh (LEFT) and the nonuniform mesh (RIGHT). Blue squares or
diamonds on the $t$-axis indicate the mesh points used in each case.
While dotted curves are the same between the left and right figures,
locations of the plotted function values (large dots on the curves)
corresponding to different meshes are significantly different. }
\end{center}

\vspace{2mm}

\no{\bf Example 2.1.} Consider meshes with twelve points over the
interval $[0, 1].$ A uniform mesh is $\TT_1=\{t_k=k\times
h\}_{k=0}^{11}$ with $h=1/11$ (blue diamonds in Figure 2.1), while a
nonuniform mesh can be
$\TT_2=\{t_0=0,~t_{k+1}=t_k+h_k,~k=0,1,\ldots,10\},$ where
$h_k,~k=0,1,\ldots,10,$ are calculated based on an equi-arclength
formula for $y=\sin(2\pi t);$ that is, $h_k,~k=1,2,\ldots,11,$ are
obtained by evenly dividing the arc of the curve of the function
over $[0, 1],$ and then vertically projecting the resulting arc
pieces with equal arclength to the $t$-axis (the mesh steps are
separated by blue squares in Figure 2.1). A composite trapezoidal
rule is employed to evaluate the arclength. Let $y_1=\sin(2\pi
t),~y_2=-0.5\sin(4\pi t),$ and $y_3=-0.8\sin(\pi t)$ be three
trigonometric functions with distinct frequencies. Their graphs over
$\TT_1$ (LEFT) and $\TT_2$ (RIGHT) are given in Figure 2.1 as red,
green and cyan dots, respectively. The true functions (continuous
curves) are shown for comparisons. {\sc Matlab$^{\circledR}$} is
used \cite{matlab}.

\vspace{3mm}

\no{\bf Definition 2.2.} Let function $y=f(t)$ be defined on the
interval $[a, b]$ in addition to $n\geq 1.$ We say that $f$ is
{\em $n$ times continuously differentiable on\/} $[a, b]$ if in
additional to the continuity of the $n$th derivative
$f^{(n)}(t),~t\in (a, b),$ both directional derivatives,
$f^{(n)}(a^+),~f^{(n)}(b^-),$ exist. We further say that $f$ is
{\em sufficiently smooth on\/} $[a, b]$ if $n$ can be as large as
we wish.

\vspace{3mm}

\no{\bf Definition 2.3.} Let functions $y=f(t)$ and $y=g(t)$ be
defined on the mesh $\TT,$ where $g$ is considered as an
approximation of $f.$ If $$|f(t)-g(t)|=O\l(h^p\r),~~~t\in \TT,$$
where $h=\max_{0\leq k\leq m}h_k,$ then we say that the
approximation at $t$ is {\em accurate to the order $p$\/} with
respect to the step sizes. From the approximation point of view,
an approximation $g$ is {\em consistent\/} if and only if $p>0~$
\cite{atkin,jain1}.

\vspace{3mm}

\no{\bf Definition 2.4.} Let $f=\{f_0,~f_1,\ldots,
f_{m+1}\},~g=\{g_0,~g_1,\ldots, g_{m+1}\}$ be two functions defined
on the mesh $\TT,$ where $g$ is considered as an approximation of
$f.$ We define the {\em scaled local difference\/} between the two
functions at $t_k$ as
$$\mbox{sld}(f,g)_k=\frac{f_k-g_k}{a},~~~
k=0,~1,\ldots,m+1,$$ where $a=\max_{0\leq k\leq m+1}|f_k|> 0.$ The
{\em scaled global error indicator\/} is defined as
$$\mbox{sgei}(f,g)=\max_{0\leq k\leq
m+1}\l|\mbox{sld}(f,g)_k\r|.$$ Note that the above definitions are
different from standard definitions of local and global relative
errors, where signs are rarely considered. The value of sld$(f,g)_k$
offers not only scaled relative error information, but also the {\em
direction of the error,\/} that is, whether $g_k$ is greater or less
than $f_k.$ The latter is particularly useful if approximations of
oscillatory problems are investigated. The function sgei$(f,g)$
provides a scaled overall error estimate and is easy to use for its
simplicity in structure. Both definitions can be used when some of
the $f_k$ values are zero.

\vspace{0.25in}

{\large\bf 3. First order
differences.~}\setcounter{section}{3}\setcounter{equation}{0} We
assume that the function $y=f(t)$ is $n$ times continuously
differentiable on $[a, b].$ Let $\TT$ be a mesh over $[a, b].$ We
define the forward, backward, and central differences as follows:
\bbb D_{+}
f(t_k)&=&\frac{f(t_{k+1})-f(t_k)}{t_{k+1}-t_k},~~~k=0,1,\ldots,m,\label{fd}\\
D_{-}
f(t_k)&=&\frac{f(t_k)-f(t_{k-1})}{t_k-t_{k-1}},~~~k=1,2,\ldots,m+1,\label{bd}\\
\delta f(t_k)&=&\frac{f(t_{k+1})-f(t_{k-1})}{t_{k+1}-t_{k-1}},~~~
k=1,2,\ldots,m. \label{cd}\eee

\vspace{3mm}

\no{\bf Theorem 3.1.} {\em Let $\TT$ be nonuniform and $f$ be twice
continuously differentiable. Then the forward, backward, and central
differences are first order approximations of the derivative
function $f'$ on $\TT.$ Further,
 \bbb D_+
f(t_k)-f'(t_k)&=&
\frac{h_k}{2}f''(\xi_1),\label{x1}\\
D_- f(t_k)-f'(t_k)&=& -\frac{h_{k-1}}{2}f''(\zeta_1),\label{y1}\\
\delta f(t_k)-f'(t_k)&=&
\frac{1}{2(h_k+h_{k-1})}\l[h_k^2f''(\xi_2)-h_{k-1}^2
f''(\zeta_2)\r],\label{z1}\eee where
$t_{k+1}>\xi_{\ell}>t_k,~t_k>\zeta_{\ell}>t_{k-1},~\ell=1,2.$ }

\vspace{2mm}

\no{\em Proof.\/} Since the proofs of the forward and backward
difference approximations are similar, we only need to show the
cases involving forward and central differences. First, according to
the Maclaurin series expansion \cite{stewart1}, we have \bbb
f(t_{k+1})&=& f(t_k)+h_kf'(t_k)+\frac{h_k^2}{2}f''(t_k)+\cdots
+\frac{h_k^n}{n!}f^{(n)}(\xi),\label{t1}\\
f(t_{k-1})&=&
f(t_k)-h_{k-1}f'(t_k)+\frac{h_{k-1}^2}{2}f''(t_k)-\cdots
+(-1)^n\frac{h_{k-1}^n}{n!}f^{(n)}(\zeta).~~~\label{t2}\eee
Therefore, letting $n=2,$\bbbb
f(t_{k+1})-f(t_k)&=&h_kf'(t_k)+\frac{h_k^2}{2}f''(\xi),\\
f(t_{k+1})-f(t_{k-1})&=&(h_k+h_{k-1})f'(t_k)+
\frac{1}{2}\l[h_k^2f''(\xi)-h_{k-1}^2f''(\zeta)\r].\eeee Equations
\R{x1} and \R{z1} become obvious. Since $h_k/(h_k+h_{k-1})<1$ and
$h_{k-1}/(h_k+h_{k-1})<1,$ according to Definition 2.3, the
differences involved are indeed first order approximations. \qed

\vspace{3mm}

\no{\bf Theorem 3.2.} {\em Let $\TT$ be uniform. Then the central
difference becomes a second order approximation of the derivative
function $f'$ if $f$ is three times continuously differentiable
while the forward and backward differences remain as first order
approximations of $f'$ if $f$ is at least twice continuously
differentiable. Further, \bbb D_+ f(t_k)-f'(t_k)&=&
\frac{h}{2}f''(\xi_1),\label{x11}\\
D_- f(t_k)-f'(t_k)&=& -\frac{h}{2}f''(\zeta_1),\label{y11}\\
\delta f(t_k)-f'(t_k)&=&
\frac{h^2}{6}\l[f'''(\xi_2)+f'''(\zeta_2)\r],\label{z11}\eee where
$t_{k+1}>\xi_{\ell}>t_k,~t_k>\zeta_{\ell}>t_{k-1},~\ell=1,2.$ }

\vspace{2mm}

\no{\em Proof.\/} The corollary is a direct extension of the Theorem
3.1 due to $h_k=h_{k-1}=h,$ Definition 2.3, and \R{t1}, \R{t2}. \qed

\vspace{3mm}

\no{\bf Remark 3.1.} Though $D_+ f(t_k)=D_- f(t_{k+1})$ for
$k=0,1,\ldots,m,$ on any mesh considered, the central difference is
not a simple combination of the forward and backward differences on
a nonuniform mesh. The central difference can be viewed as an
arithmetic average of the forward and backward differences only on
any uniform mesh. Interestingly, this average improves the accuracy
of the approximation significantly.

\vspace{3mm}

\no{\bf Remark 3.2.} The right-hand-side of \R{x1}-\R{z1} and
\R{x11}-\R{z11} can be viewed as errors of the respective
approximations.

\vspace{0.25in}

{\large\bf 4. Second order
differences.~}\setcounter{section}{4}\setcounter{equation}{0} Let
function $y=f(t)$ be $n$ times continuously differentiable on $[a,
b],$ where $n$ is sufficiently large. We will investigate if we
can approximate the second order derivative function by applying
the forward, backward and central differences repeatedly.

\vspace{3mm}

\no{\bf Theorem 4.1.} {\em Let $\TT$ be nonuniform. Then
$$P(Qf(t_k))\neq Q(Pf(t_k))~\mbox{ for any applicable $t_k\in\TT,$}$$
where $P$ and $Q$ are any two different difference operations
denoted by $D_+,~D_-$ or $\delta.$ }

\vspace{2mm}

\no{\em Proof.\/} We only need to show the cases involving $D_+(D_-
f(t_k))$ and $D_+(\delta f(t_k)).$ According to \R{fd} and \R{bd},
we have \bbb D_+(D_- f(t_k))&=&\frac{D_- f(t_{k+1})-D_-
f(t_k)}{h_k}=\l[\frac{f(t_{k+1})-f(t_k)}{h_k}-\frac{f(t_k)-f(t_{k-1})}{h_{k-1}}\r]/{h_k}\nnn\\
&=&\frac{h_{k-1}f(t_{k+1})-(h_k+h_{k-1})f(t_k)+h_kf(t_{k-1})}{h_k^2h_{k-1}}.\label{fb}\eee
By the same token, \bbb D_-(D_+ f(t_k))&=&\frac{D_+ f(t_k)-D_+
f(t_{k-1})}{h_{k-1}}=\l[\frac{f(t_{k+1})-f(t_k)}{h_k}-\frac{f(t_k)-f(t_{k-1})}{h_{k-1}}\r]/{h_{k-1}}\nnn\\
&=&\frac{h_{k-1}f(t_{k+1})-(h_k+h_{k-1})f(t_k)+h_kf(t_{k-1})}{h_kh_{k-1}^2}.\label{bf}\eee
Recalling that $h_k\neq h_{k-1},$ we have
$$D_+(D_- f(t_k))\neq D_-(D_+ f(t_k)).$$
Further, \bbb &&D_+(\delta f(t_k))=\frac{\delta f(t_{k+1})-\delta
f(t_k)}{h_k}=\l[\frac{f(t_{k+2})-f(t_k)}{h_{k+1}+h_k}-
\frac{f(t_{k+1})-f(t_{k-1})}{h_k+h_{k-1}}\r]/{h_k}\nnn\\
&&~~~~~~=\frac{(h_k+h_{k-1})(f(t_{k+2})-f(t_k))-
(h_{k+1}+h_k)(f(t_{k+1})-f(t_{k-1}))}{(h_{k+1}+h_k)h_k(h_k+h_{k-1})}.\label{fc}\\
&&\delta(D_+ f(t_k))=\frac{D_+ f(t_{k+1})-D_+
f(t_{k-1})}{h_k+h_{k-1}}\nnn\\&&~~~~~~=
\l[\frac{f(t_{k+2})-f(t_{k+1})}{h_{k+1}}-\frac{f(t_k)-f(t_{k-1})}{h_{k-1}}\r]/{(h_k+h_{k-1})}\nnn\\
&&~~~~~~=\frac{h_{k-1}(f(t_{k+2})-f(t_{k+1}))-
h_{k+1}(f(t_k)-f(t_{k-1}))}{h_{k+1}(h_{k-1}+h_k)h_{k-1}}.\label{cf}\eee
Therefore
$$\delta(D_+ f(t_k))\neq D_+(\delta f(t_k))$$
unless $h_{k-1}=h_k=h_{k+1}.$ \qed

\vspace{3mm}

\no{\bf Corollary 4.1.} {\em Let $\TT$ be uniform. Then we have
$$P(Qf(t_k))= Q(Pf(t_k))~\mbox{ for all applicable $t_k\in\TT,$}$$
where $P$ and $Q$ are any two different difference operations
denoted by $D_+,~D_-$ or $\delta.$ }

\vspace{2mm}

The proof of the corollary follows from \R{fb}-\R{cf}. To explore
interesting features of the proposed formulas on different meshes,
we require exceptionally high orders of the derivative functions in
the following theorem. The requirements can be conveniently eased
when the types of the meshes are fixed.

\vspace{3mm}

\no{\bf Theorem 4.2.} {\em Let $\TT$ be nonuniform. Then none of the
second order differences $P(Qf(t_k)),$ where $P$ and $Q$ are any two
difference operations denoted by $D_+,~D_-$ or $\delta,$ is a
consistent approximation of $f''(t_k).$ Further, if $f$ is
sufficiently smooth, \bbb D_+(D_+ f(t_k))&=&
\frac{h_{k+1}+h_k}{2h_k}f''(t_k)+\frac{
(h_{k+1}+h_k)(h_{k+1}+2h_k)}{6h_k}f'''(t_k) \nnn\\&&+
\frac{h_{k+1}+h_k}{24h_kh_{k+1}}\l[(h_{k+1}+h_k)^3f^{(4)}(\tilde{\xi})
-h_k^3f^{(4)}(\xi)\r],\label{xx2}\\ D_-(D_- f(t_k))&=&
\frac{h_{k-1}+h_{k-2}}{2h_{k-1}}f''(t_k)
-\frac{(h_{k-1}+h_{k-2})(2h_{k-1}+h_{k-2})}{6h_{k-1}}f'''(t_k)\nnn\\
&&
+\frac{h_{k-1}+h_{k-2}}{24h_{k-1}h_{k-2}}\l[-h_{k-1}^3f^{(4)}(\zeta)+(h_{k-1}+h_{k-2})^3
f^{(4)}(\tilde{\zeta})\r],\label{yy2}\\
\delta(\delta f(t_k))&=&
\frac{h_{k+1}+h_k+h_{k-1}+h_{k-2}}{2(h_k+h_{k-1})}f''(t_k)\nnn\\&&+
\frac{(h_{k+1}+h_k)^2-(h_{k-1}+h_{k-2})^2}{6(h_k+h_{k-1})}f'''(t_k)
\nnn\\
&&+\frac{(h_k+h_{k+1})^3+(h_{k-1}+h_{k-2})^3}{24(h_k+h_{k-1})}f^{(4)}(t_k)\nnn\\
&&+\frac{(h_{k+1}+h_k)^4f^{(5)}(\tilde{\xi})-
(h_{k-1}+h_{k-2})^4f^{(5)}(\tilde{\zeta})}{120(h_k+h_{k-1})}, \label{zz2}\\
D_+(D_- f(t_k))&=&
\frac{h_k+h_{k-1}}{2h_k}f''(t_k)+\frac{h_k^2-h_{k-1}^2}{6h_k}f'''(t_k)
+\frac{h_k^3+h_{k-1}^3}{24h_k}f^{(4)}(t_k)\nnn\\&&+
\frac{1}{120h_k}\l[h_k^4f^{(5)}(\xi)-h_{k-1}^4f^{(5)}(\zeta)\r],\label{xy2}\\
D_-(D_+
f(t_k))&=&\frac{h_k+h_{k-1}}{2h_{k-1}}f''(t_k)+\frac{h_k^2-h_{k-1}^2}{6h_{k-1}}f'''(t_k)
+\frac{h_k^3+h_{k-1}^3}{24h_{k-1}}f^{(4)}(t_k)\nnn\\&&~~~~~~+\frac{1}{120h_{k-1}}\l[h_k^4f^{(5)}(\xi)
-h_{k-1}^4f^{(5)}(\zeta)\r], \label{yx2}\\ D_+(\delta f(t_k))&=&
\frac{h_{k+1}+h_{k-1}}{2h_k}f''(t_k)+\frac{(h_{k+1}+h_k)^2-h_k^2+h_kh_{k-1}
-h_{k-1}^2}{6h_k}f'''(t_k)\nnn\\
&&+\frac{1}{24h_k(h_k+h_{k-1})}\l[
(h_k+h_{k-1})(h_{k+1}+h_k)^3f^{(4)}(\tilde{\xi})\r.\nnn\\&&~~~~~~\l.-h_k^4f^{(4)}(\xi)
-h_{k-1}^4f^{(4)}(\zeta)\r],\label{xz2}\\ \delta(D_+ f(t_k))&=&
\frac{h_{k+1}+2h_k+h_{k-1}}{2(h_k+h_{k-1})}f''(t_k)+\frac{(h_k+h_{k+1})^3-h_k^3-h_{k+1}h_{k-1}^2}{6
h_{k+1}(h_k+h_{k-1})}f'''(t_k)\nnn\eee\bbb &&
+\frac{1}{24h_{k+1}(h_k+h_{k-1})}\l[(h_{k+1}+h_k)^4f^{(4)}(\tilde{\xi})-h_k^4f^{(4)}(\xi)\r.\nnn\\
&&\l.h_{k-1}^3h_{k+1}f^{(4)}(\zeta)\r],\label{zx2}\eee\bbb
D_-(\delta f(t_k))&=&
\frac{h_k+h_{k-2}}{2h_{k-1}}f''(t_k)+\frac{h_k^3+h_{k-1}^3-(h_k+h_{k-1})(h_{k-1}+h_{k-2})^2}{6(h_k+h_{k-1})
h_{k-1}}f'''(t_k)\nnn\\
&&+\frac{1}{24(h_k+h_{k-1})
h_{k-1}}\l[h_k^4f^{(4)}(\xi)+(h_{k-1}+h_{k-2})^4f^{(4)}(\zeta)\r.\nnn\\&&\l.+(h_k+h_{k-1})
(h_{k-1}+h_{k-2})^3f^{(4)}(\tilde{\zeta})\r], \label{yz2}\\
\delta(D_-
f(t_k))&=&\frac{h_k+2h_{k-1}+h_{k-2}}{2(h_k+h_{k-1})}f''(t_k)+
\frac{h_{k-2}h_k^2+h_{k-1}^3-(h_{k-1}+h_{k-2})^3}{6h_{k-2}(h_k+h_{k-1})}f'''(t_k)\nnn\\
&&
+\frac{1}{24(h_k+h_{k-1})h_{k-2}}\l[h_{k-2}h_k^3f^{(4)}(\xi)-h_{k-1}^4f^{(4)}(\zeta)\r.\nnn\\
&&\l.+(h_{k-1}+h_{k-2})^4f^{(4)}(\tilde{\zeta})\r],\label{zy2}\eee
where $\tilde{\xi},~\xi,~\tilde{\zeta}$ and $\zeta$ are different
numbers and
$t_{k+2}>\tilde{\xi}>t_k,~t_{k+1}>\xi>t_k,~t_k>\zeta>t_{k-1},~t_k>\tilde{\zeta}>t_{k-2}.$
 }

\vspace{2mm}

\no{\em Proof.\/} Due to their similarities, here we only present
proofs of a few cases as an illustration. Let us consider $D_+(D_+
f(t_k)),~\delta(\delta f(t_k)),$ and $D_-(\delta f(t_k)).$ It will
be sufficient to show \R{xx2}, \R{zz2} and \R{yz2}. We have the
Maclaurin series \bbb f(t_{k+2})&=&
f(t_k)+(h_{k+1}+h_k)f'(t_k)+\frac{(h_{k+1}+h_k)^2}{2}f''(t_k)
\nnn\\
&&
+\frac{(h_{k+1}+h_k)^3}{6}f'''(t_k)+\cdots+\frac{(h_{k+1}+h_k)^n}{n!}f^{(n)}(\tilde{\xi}),\label{t3}\\
f(t_{k-2})&=&
f(t_k)-(h_{k-1}+h_{k-2})f'(t_k)+\frac{(h_{k-1}+h_{k-2})^2}{2}f''(t_k)\nnn\\
&&-
\frac{(h_{k-1}+h_{k-2})^3}{6}f'''(t_k)+\cdots+(-1)^n\frac{(h_{k-1}+h_{k-2})^n}{n!}f^{(n)}(\tilde{\zeta}).~~~~~~\label{t4}\eee
Similar to \R{fb} and \R{bf}, we have \bb{ff}D_+(D_+
f(t_k))=\frac{h_kf(t_{k+2})-(h_{k+1}+h_k)f(t_{k+1})+h_{k+1}f(t_k)}{h_{k+1}h_k^2}.\ee
Now, according to \R{t1}, \R{t3}, \bbbb
&&\hspace{-6mm}h_kf(t_{k+2})-(h_{k+1}+h_k)f(t_{k+1})+h_{k+1}f(t_k)\\&&\hspace{-6mm}~~~=h_k\l(f(t_k)
+(h_{k+1}+h_k)f'(t_k)+\frac{(h_{k+1}+h_k)^2}{2}f''(t_k)+\cdots
+\frac{(h_{k+1}+h_k)^n}{n!}f^{(n)}(\tilde{\xi})\r)\\
&&\hspace{-6mm}~~~~~-(h_{k+1}+h_k)\l(f(t_k)
+h_kf'(t_k)+\frac{h_k^2}{2}f''(t_k)+\cdots+\frac{h_k^n}{n!}f^{(n)}(\xi)\r)+h_{k+1}f(t_k)\eeee\bbbb
&&\hspace{-6mm}~~~=\frac{h_kh_{k+1}(h_{k+1}+h_k)}{2}f''(t_k)+\frac{h_k(h_{k+1}+h_k)
\l[(h_{k+1}+h_k)^2-h_k^2\r]}{3!}f'''(t_k)
+\cdots\\&&\hspace{-6mm}~~~~~+
\frac{h_k(h_{k+1}+h_k)}{n!}\l[(h_{k+1}+h_k)^{n-1}f^{(n)}(\tilde{\xi})-h_k^{n-1}f^{(n)}(\xi)\r]\\
&&\hspace{-6mm}~~~=\frac{h_kh_{k+1}(h_{k+1}+h_k)}{2}f''(t_k)+\frac{h_kh_{k+1}(h_{k+1}+h_k)
(h_{k+1}+2h_k)}{3!}f'''(t_k) +\cdots\\&& \hspace{-6mm}~~~~~+
\frac{h_k(h_{k+1}+h_k)}{n!}\l[(h_{k+1}+h_k)^{n-1}f^{(n)}(\tilde{\xi})-h_k^{n-1}f^{(n)}(\xi)\r].\eeee
Substituting the above into \R{ff} and letting $n=4,$ we obtain
\R{xx2} which indicates that the difference is not a consistent
approximation of $f''(t_k)$ unless $h_{k+1}=h_k.$ Our next case
involving \R{zz2} is interesting. To see this, we have \bbbb
&&\hspace{-6mm}\delta(\delta f(t_k))=\frac{\delta f(t_{k+1})-\delta
f(t_{k-1})}{h_k+h_{k-1}}\nnn\\&&\hspace{-6mm}~~=
\l[\frac{f(t_{k+2})-f(t_k)}{h_{k+1}+h_k}-\frac{f(t_k)-f(t_{k-2})}{h_{k-1}+h_{k-2}}\r]/{(h_k+h_{k-1})}\nnn\\
&&\hspace{-6mm}~~=\frac{(h_{k-1}+h_{k-2})f(t_{k+2})-(h_{k+1}+h_k+h_{k-1}+h_{k-2})f(t_k)+
(h_{k+1}+h_k)f(t_{k-2})}{(h_{k+1}+h_k)(h_k+h_{k-1})(h_{k-1}+h_{k-2})}.\eeee
We obtain \R{zz2} readily from the above by substituting \R{t3} and
\R{t4}. Further, using steps similar to those in the above
discussion, we have \bbbb D_-(\delta
f(t_k))&=&\l[(h_{k-1}+h_{k-2})f(t_{k+1})-
(h_k+h_{k-1})f(t_k)-(h_{k-1}+h_{k-2})f(t_{k-1})\r.\\
&&\l.+(h_k+h_{k-1})f(t_{k-2})\r]/[(h_k+h_{k-1})h_{k-1}(h_{k-1}+h_{k-2})].\eeee
Substitute \R{t1}, \R{t2} and \R{t4} into the above equation. We
acquire \R{yz2} after simplification. Obviously, \R{yz2} indicates
again that the finite difference is not a consistent approximation
of $f''(t_k)$ unless $h_k=h_{k-1}=h_{k-2}.$ \qed

\vspace{3mm}

\no{\bf Corollary 4.2.} {\em Let $\TT$ be any mesh. Then
$$D_+(D_- f(t_k))=\frac{h_{k-1}}{h_k}D_-(D_+ f(t_k)),~~~t_k\in\TT.$$
}

\vspace{3mm}

\no{\bf Corollary 4.3.} {\em Let $\alpha$ be a positive constant,
and let $h_{\ell+1}=\alpha h_{\ell}$ for all possible $\ell$ on a
nonuniform mesh $\TT.$ Then a necessary and sufficient condition for
any aforementioned finite difference formula to be a consistent
approximation of $f''$ is $$\alpha\equiv 1.$$ }

\vspace{1mm}

\no{\em Proof.\/} Because the proofs of different cases are similar,
we will show only one of them. Consider $\delta(D_- f(t_k)).$ To
have a consistent approximation, we must require the coefficient of
$f''(t_k)$ in \R{zy2} to be one, that is,
$$\frac{h_k+2h_{k-1}+h_{k-2}}{2(h_k+h_{k-1})}
=\frac{\alpha^2h_{k-2}+2\alpha h_{k-2}+h_{k-2}}{2(\alpha^2
h_{k-2}+\alpha h_{k-2})}=\frac{\alpha^2+2\alpha +1}{2(\alpha^2
+\alpha)}=1.$$ Solving the above equation, we acquire
$\alpha^2+2\alpha +1=2\alpha^2 +2 \alpha$ which implies that
$\alpha^2=1.$ Since $\alpha$ is positive, therefore $\alpha\equiv 1$
is our only solution. \qed

\vspace{3mm}

\no{\bf Theorem 4.3.} {\em Let $\TT$ be uniform. Then for any valid
index $k,$ $\delta(\delta f(t_k)),~D_+(D_- f(t_k))$ and $D_-(D_+
f(t_k))$ are second order approximations of $f''(t_k)$ if $f$ is
four times continuously differentiable, and all other second order
differences are first order approximations of $f''(t_k)$ if $f$ is
three times continuously differentiable. Further, \bbbb
\delta(\delta f(t_k))- f''(t_k)&=&
\frac{h^2}{6}\l[f^{(4)}(\tilde{\xi})+f^{(4)}(\tilde{\zeta})\r],
\\ D_+(D_- f(t_k))- f''(t_k)&=&
\frac{h^2}{24}\l[f^{(4)}(\xi)+f^{(4)}(\zeta)\r],\\ D_-(D_+
f(t_k))-f''(t_k)&=&D_+(D_- f(t_k))- f''(t_k),\\
D_+(D_+ f(t_k))- f''(t_k)&=& \frac{h}{3}\l[4f'''(\tilde{\xi}) -f'''(\xi)\r],\\
D_-(D_- f(t_k))- f''(t_k)&=& \frac{h}{3}\l[
f'''(\tilde{\zeta})-4f'''(\zeta)\r],\\ D_+(\delta f(t_k))-
f''(t_k)&=&\frac{h}{12}\l[ 8f'''(\tilde{\xi})-f'''(\xi)
-f'''(\zeta)\r],\\
\delta(D_+ f(t_k))- f''(t_k)&=&D_+(\delta
f(t_k))- f''(t_k),\\
D_-(\delta f(t_k))- f''(t_k)&=&
\frac{h}{12}\l[f'''(\xi)+f'''(\zeta)-8f'''(\tilde{\zeta})\r], \\
\delta(D_- f(t_k))-f''(t_k)&=& D_-(\delta f(t_k))- f''(t_k),\eeee
where $\tilde{\xi},~\xi,~\tilde{\zeta}$ and $\zeta$ are constants
for which
$t_{k+2}>\tilde{\xi}>t_k,~t_{k+1}>\xi>t_k,~t_k>\zeta>t_{k-1},~t_k>\tilde{\zeta}>t_{k-2}.$
}

\vspace{2mm}

\no{\em Proof.\/} We only need to show the nine equations listed.
This can be implemented by letting $h_{k+1}=h_k=h_{k-1}=h_{k-2}=h$
in \R{xx2}-\R{zy2} and choosing proper values of $n$ in the
Maclaurin series \R{t1}, \R{t2}, \R{t3} and \R{t4}. \qed

\vspace{3mm}

\no{\bf Remark 4.1.} The difference formula \R{aa3} is a consistent
first order approximation of $f''(t)$ on a uniform mesh.

\vspace{3mm}

\no{\bf Remark 4.2.} The right-hand sides of the nine equations in
Theorem 4.3 can be viewed as errors of the respective approximations
on uniform meshes.

\vspace{3mm}

\no{\bf Remark 4.3.} Although none of the standard second order
differences discussed in this paper may approximate $f''$ on a
nonuniform mesh, specially formulated finite difference formulas may
work. An example is \bb{good} {\cal D}_2 f(t_k)=\frac{D_+ f(t_k)-D_-
f(t_k)}{(h_{k-1}+h_k)/2}, ~~~t_k\in\TT.\ee It can be shown by using
Maclaurin series that \R{good} provides a first order approximation
of $f''(t_k)$ on any nonuniform mesh.

\vspace{0.25in}

{\large\bf 5. Numerical
examples.~}\setcounter{section}{5}\setcounter{equation}{0} It has
been known that trigonometric functions $y=\sin (at),~y=\cos (bt)$
possess excellent smoothness for applications. The functions are
infinitely differentiable. To illustrate our results, we choose the
function \bb{e1}y=-\sin(4\pi t),~~~0\leq t\leq 1.\ee The
corresponding derivatives are \bb{e2}y'=-4\pi\cos(4\pi
t),~y''=(4\pi)^2\sin(4\pi t),~~~0\leq t\leq 1.\ee

\vspace{1mm}

\begin{center}
\epsfig{file=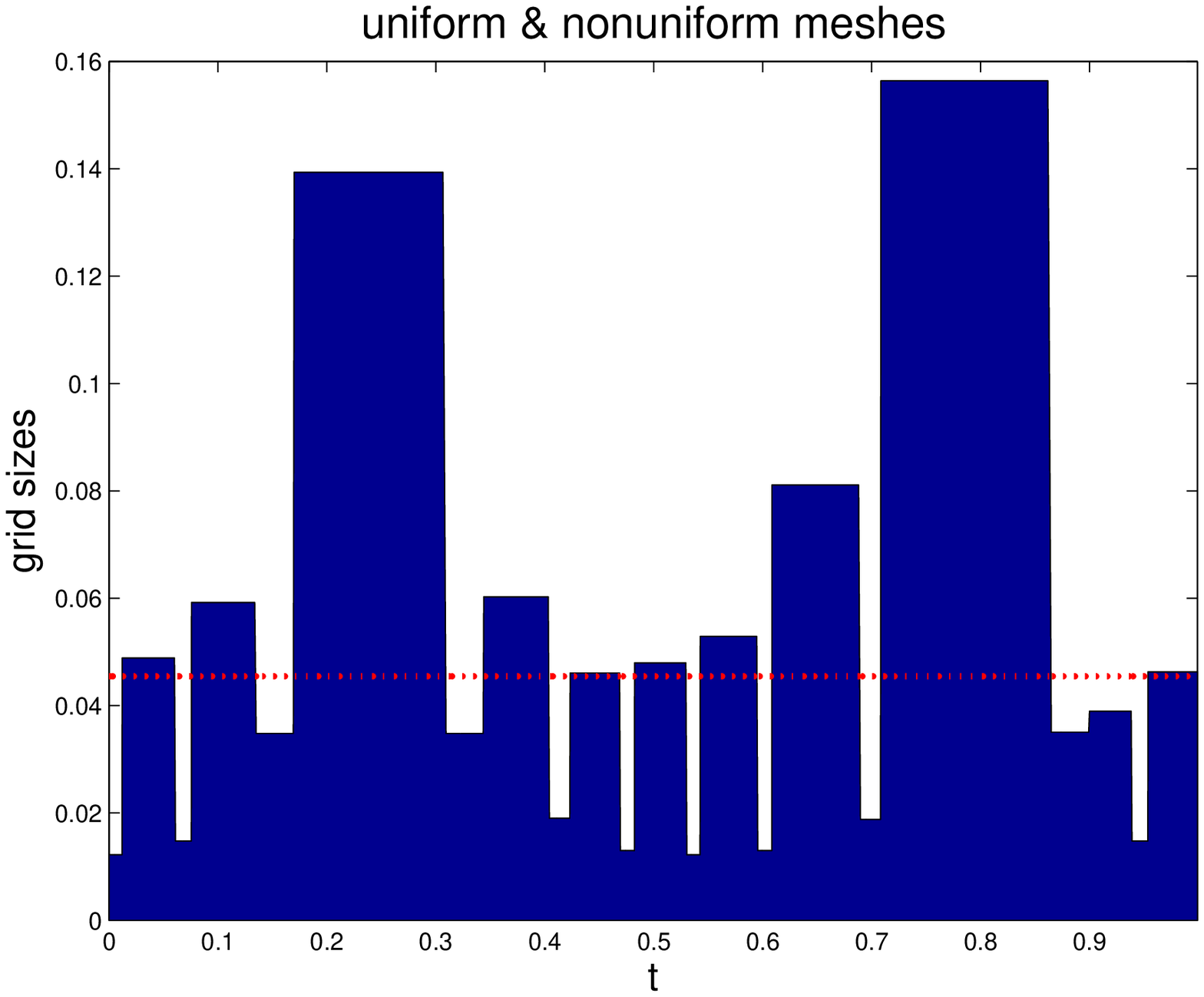,width=2.6in,height=1.6in}
\epsfig{file=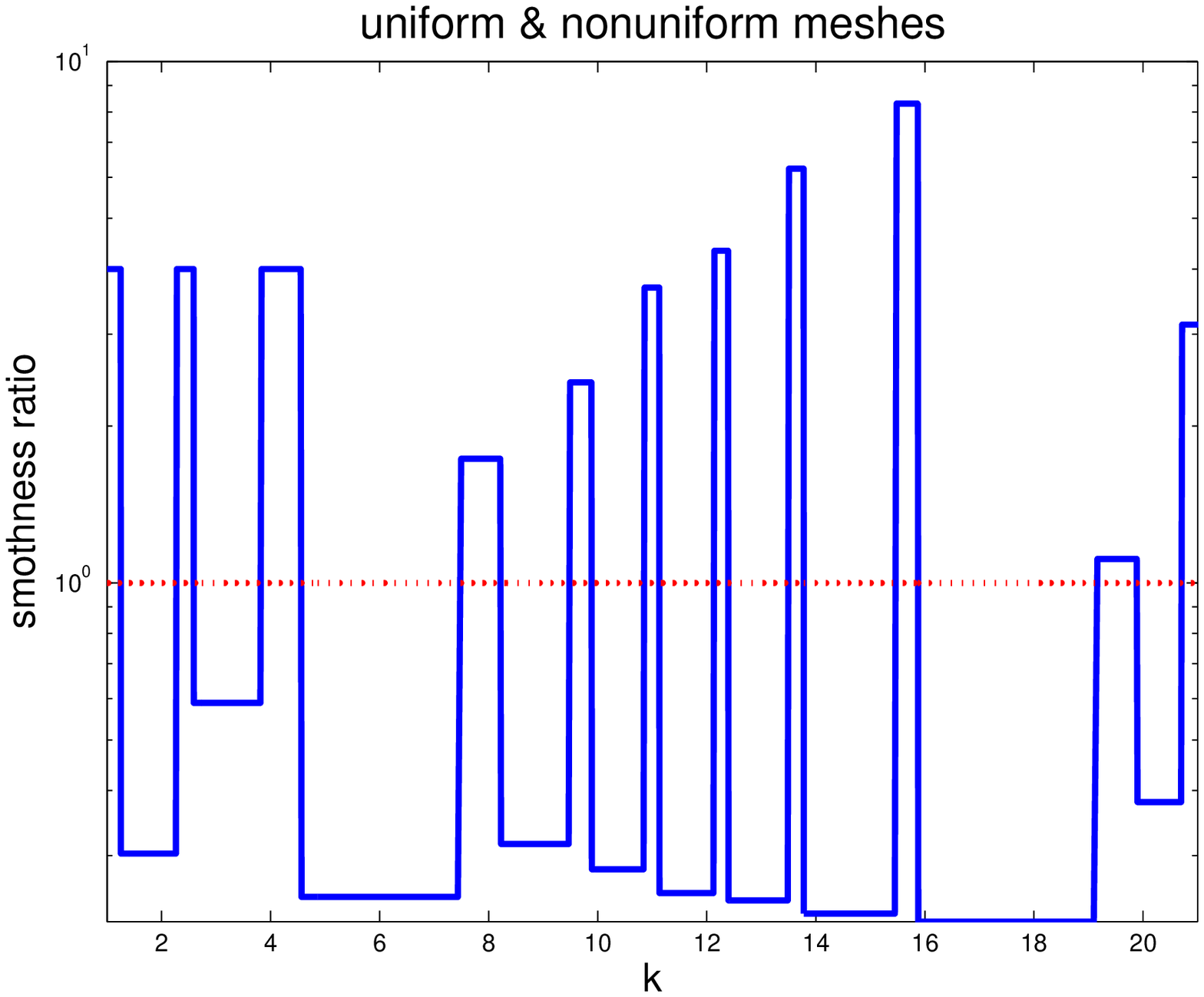,width=2.6in,height=1.6in}

\vspace{2mm}

\parbox[t]{13.5cm}{\footnotesize{\bf Figure 5.1.} Plot of
grid step sizes (red dotted line is for the uniform mesh and blue
solid curve and area are for the nonuniform mesh, LEFT). In this
figure, $[0,1]$ is divided into 22 subregions proportional to the 22
steps used. Height of each of the bars located in the subregion is
the actual size of the corresponding step size; The second figure is
for smoothness ratio for the uniform mesh (red dotted line) and the
nonuniform mesh (blue curve). In the case, $[0, 1]$ is divided into
21 subregions proportional to the first 21 steps used. Height of
each bar is for the corresponding ratio.}
\end{center}

To achieve better simulation resolution, we add eleven additional
external points into the uniform and nonuniform meshes introduced
in Example 2.1, respectively. The new points are separated by the
original points. The new points in the nonuniform mesh are
particularly set to be closer to their right-side points, as
plotted in Figure 5.1 (LEFT). In the second frame of Figure 5.1,
we show the ratio $h_{k+1}/h_k,~k=1,2,\ldots,21,$ (the red line is
for the uniform mesh and the blue curve is for the nonuniform
mesh). A logarithmic scale in the $Y$-direction is used to show
more precisely the ratio, which is often called the {\em
smoothness ratio\/} in engineering computations and is chosen
between 0 and 10.

Numerical computations will be carried out for a number of the
differences. In our experiments, a numerical solution $g$ is
acceptable if $\mbox{sgei}(f,g)<k\ll 1,$ and it is unacceptable if
$\mbox{sgei}(f,g)\geq 1.$ {\sc Matlab$^{\circledR}$} will be used
throughout the experiments.

\vspace{3mm}

\no{\bf Example 5.1.} Consider the first difference $\delta f.$ The
true first derivative in \R{e2}, and the central difference are
plotted in Figure 5.2. Though the approximations in both cases are
consistent, we may observe that the approximation in the uniform
case is much better than that in the nonuniform case, since in the
former situation the approximation is second order.

\begin{center}
\epsfig{file=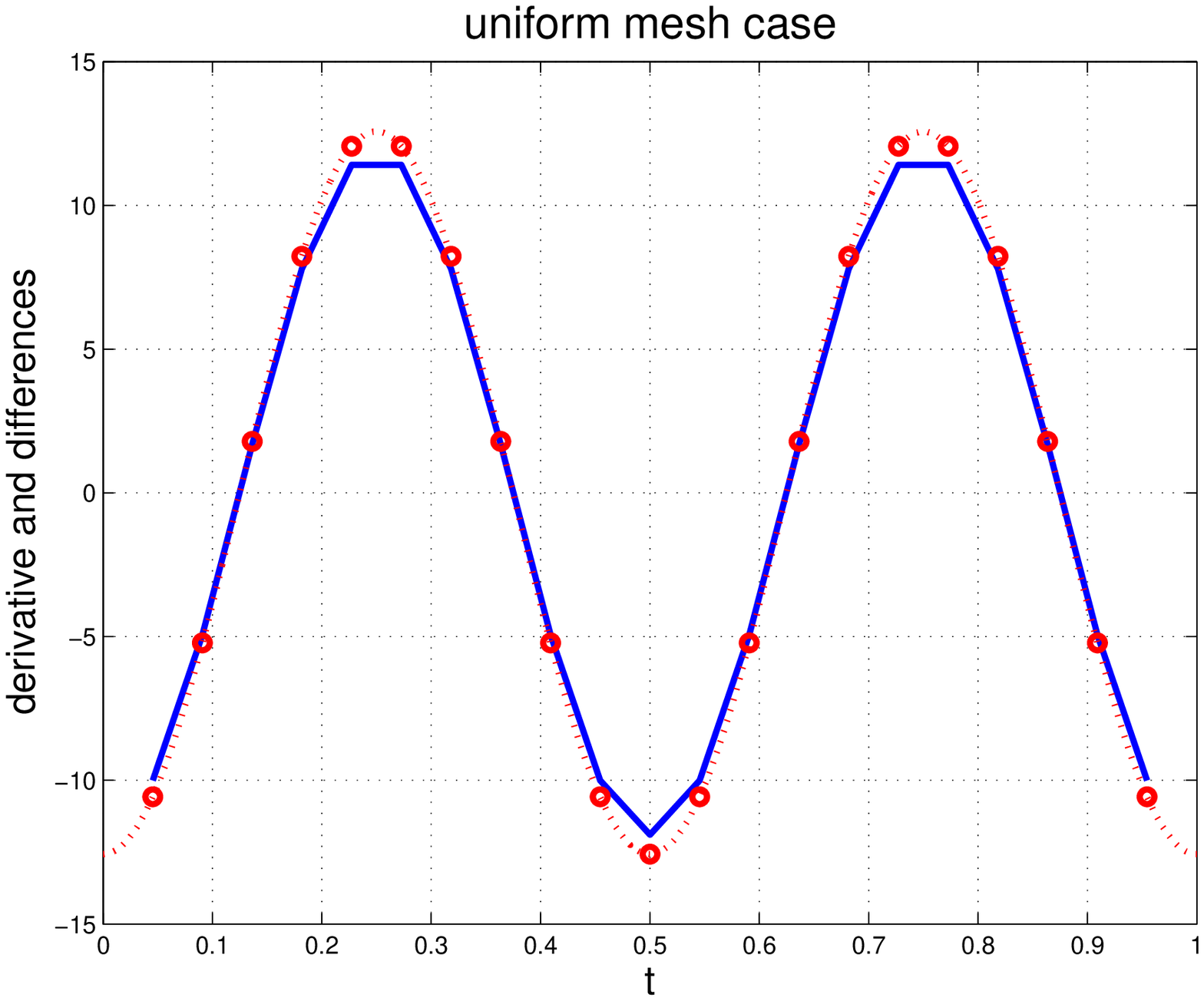,width=1.8in,height=1.6in}
\epsfig{file=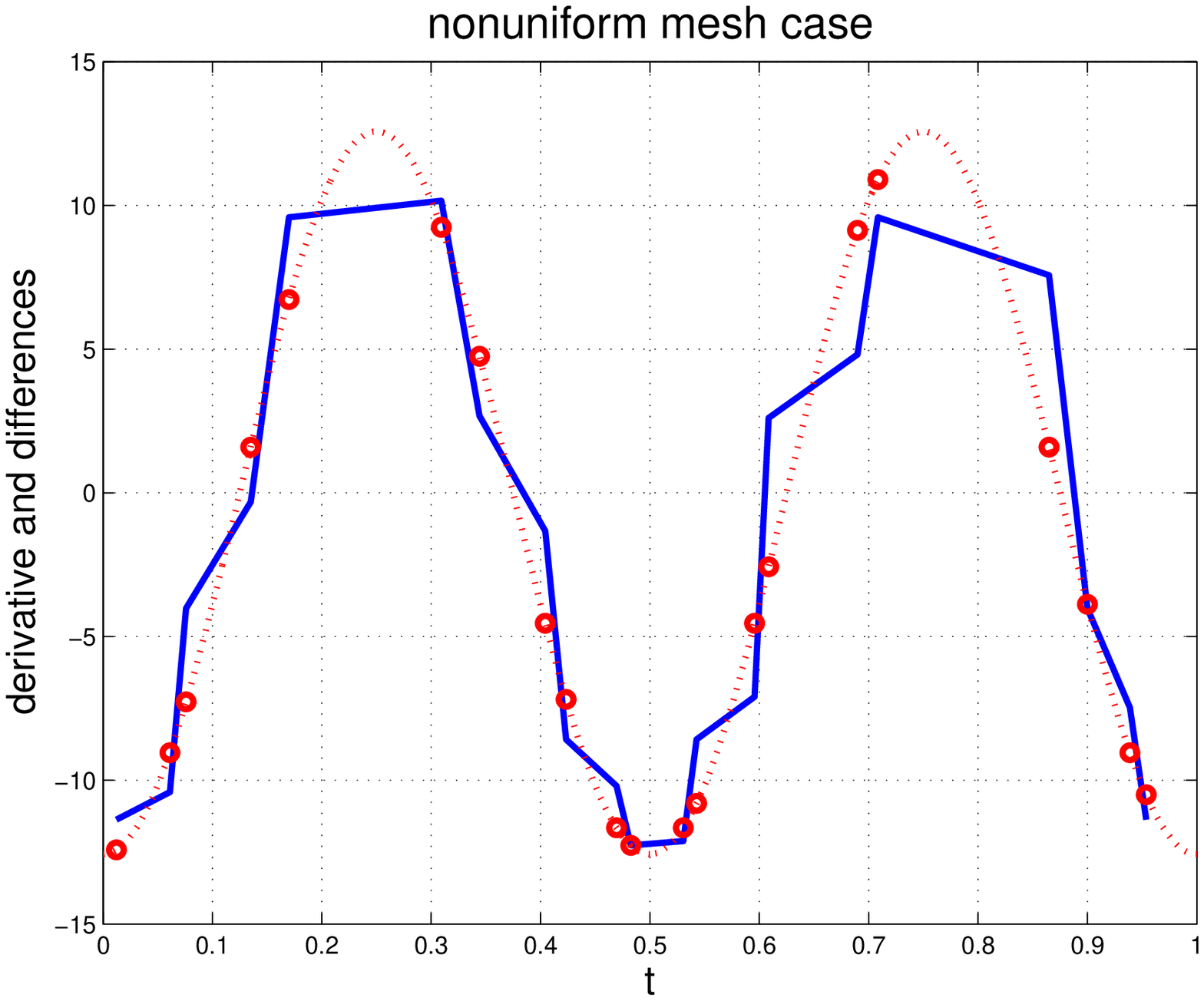,width=1.8in,height=1.6in}
\epsfig{file=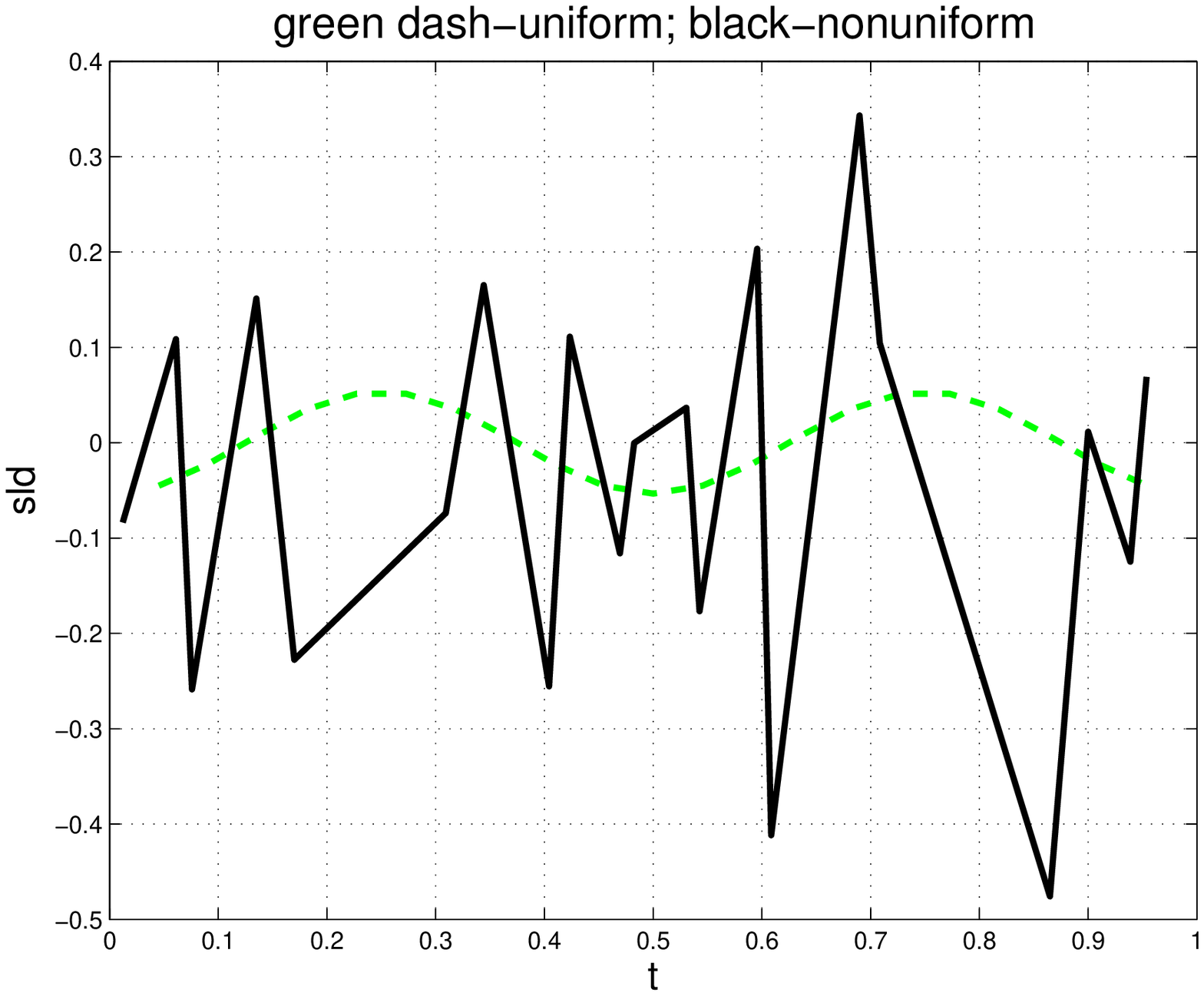,width=1.8in,height=1.6in}

\vspace{3mm}

\parbox[t]{13.5cm}{\footnotesize{\bf Figure 5.2.} Plot of
true solution (red dotted curves) and approximations (blue curves)
on the uniform mesh (LEFT) and the nonuniform mesh (MIDDLE).
Corresponding scaled local difference curves are given (RIGHT: green
for the uniform mesh case and black for the nonuniform case). We may
observe from the first two pictures that, even though approximations
on both meshes are consistent, the approximation on the uniform mesh
is much ``nicer'' than that on the nonuniform mesh. The third
picture confirms this by showing the scaled difference values over
the domain (the black curve acts more violently with a relatively
large amplitude).}
\end{center}

These are further confirmed by the third frame in Figure 5.2 in
which relative errors are given. While the scaled global error
indicator in the uniform mesh case $\mbox{sgei}\approx 0.0535,$ in
the nonuniform case, it reaches 0.4759 which is about 9 times the
uniform mesh case! The irregular smoothness ratios shown in Figure
5.1 may explain why the error in the nonuniform grid case is
oscillatory.

\vspace{3mm}

\no{\bf Example 5.2.} Consider the second difference $D_+ (D_+ f).$
We plot the true second derivative given in \R{e2} and the
difference in Figure 5.3. It can be observed that while the
difference approximates the derivative function reasonably over the
uniform mesh (global error indicator value $\mbox{sgei}\approx
0.5559$ because relatively large steps are used), it produces an
unacceptable results on the nonuniform mesh with an indicator value
$\mbox{sgei}\approx 2.5487.$ The irregular oscillations of the
finite difference on nonuniform mesh is also unacceptable.

\begin{center}
\epsfig{file=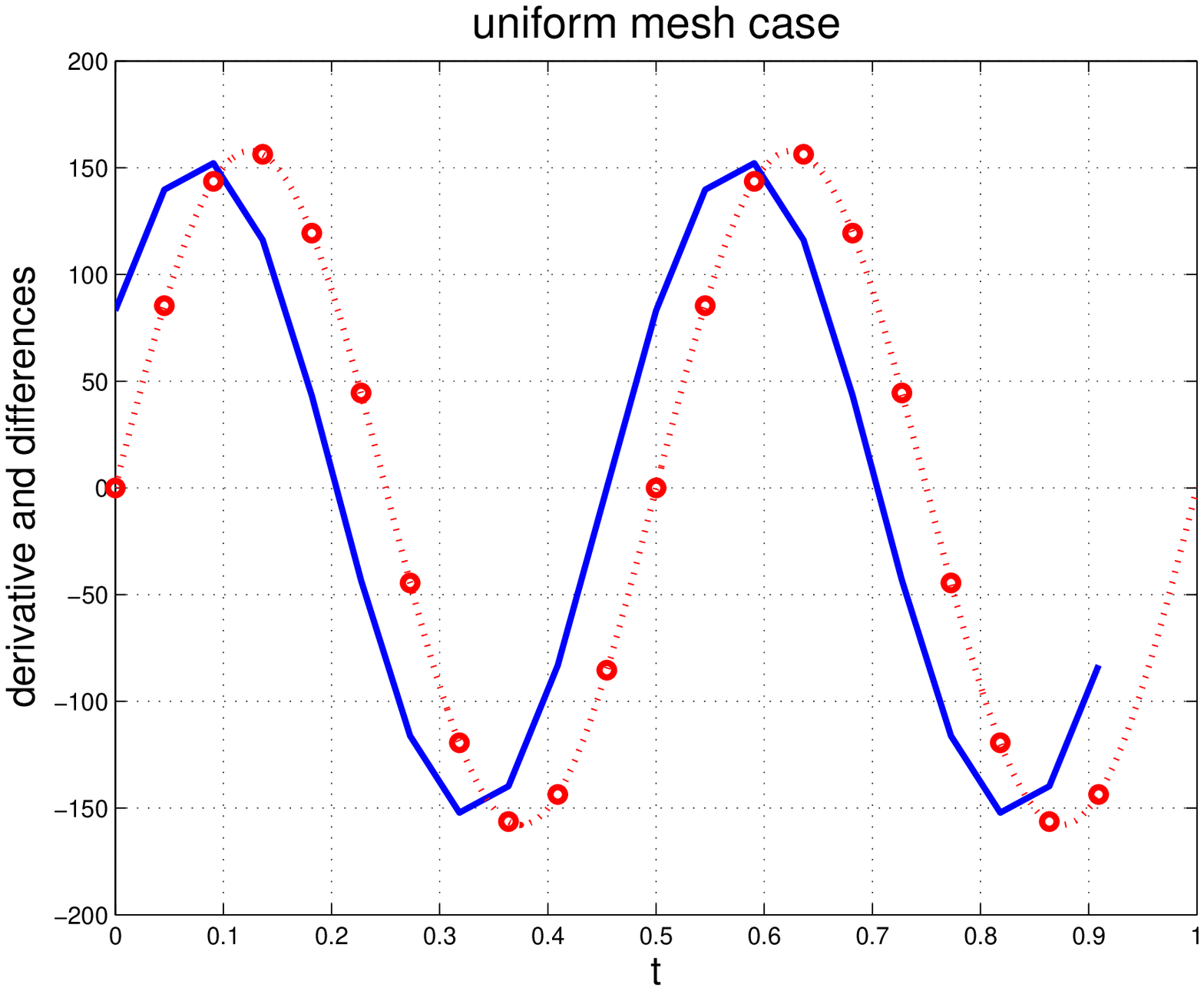,width=1.8in,height=1.6in}
\epsfig{file=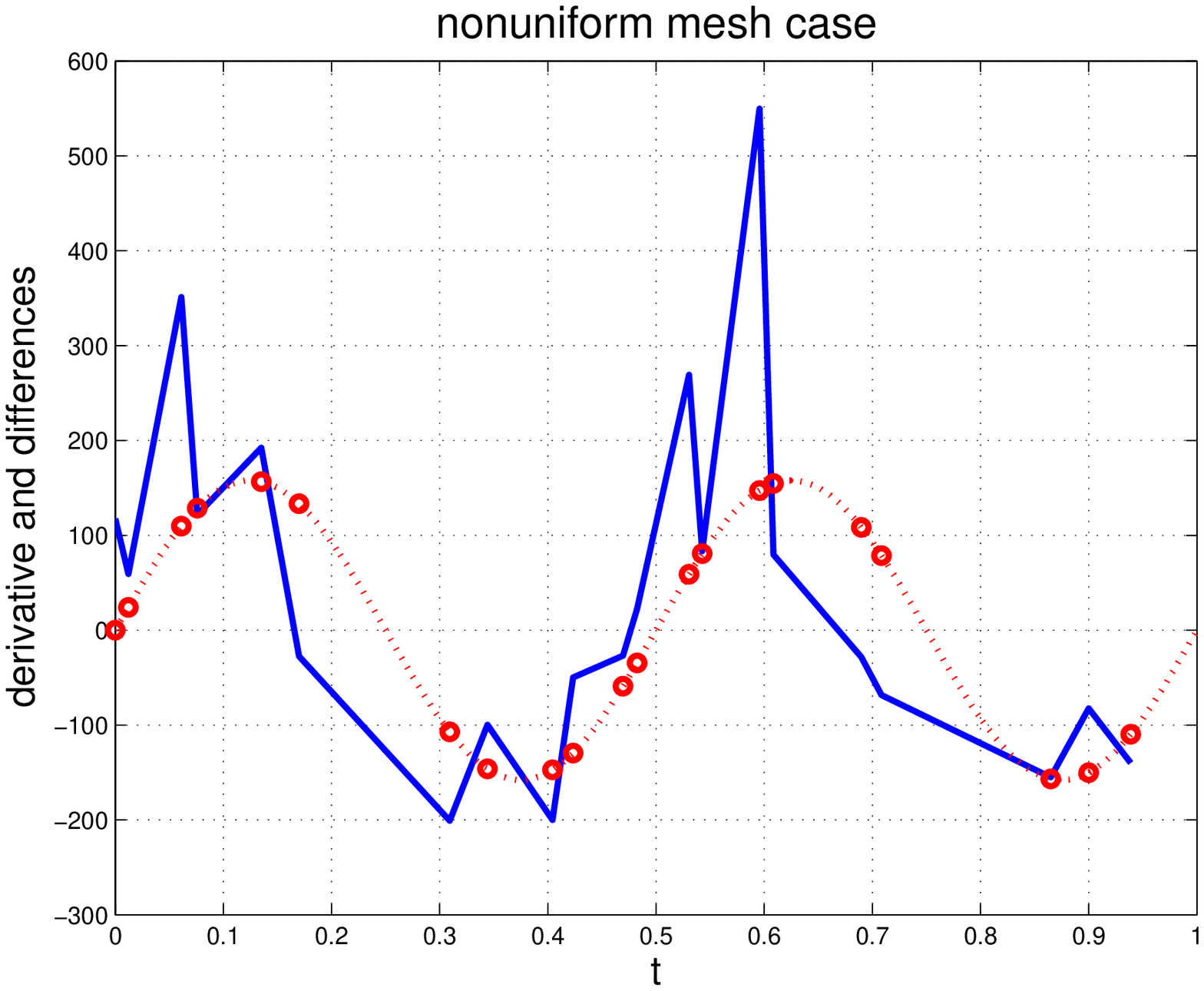,width=1.8in,height=1.6in}
\epsfig{file=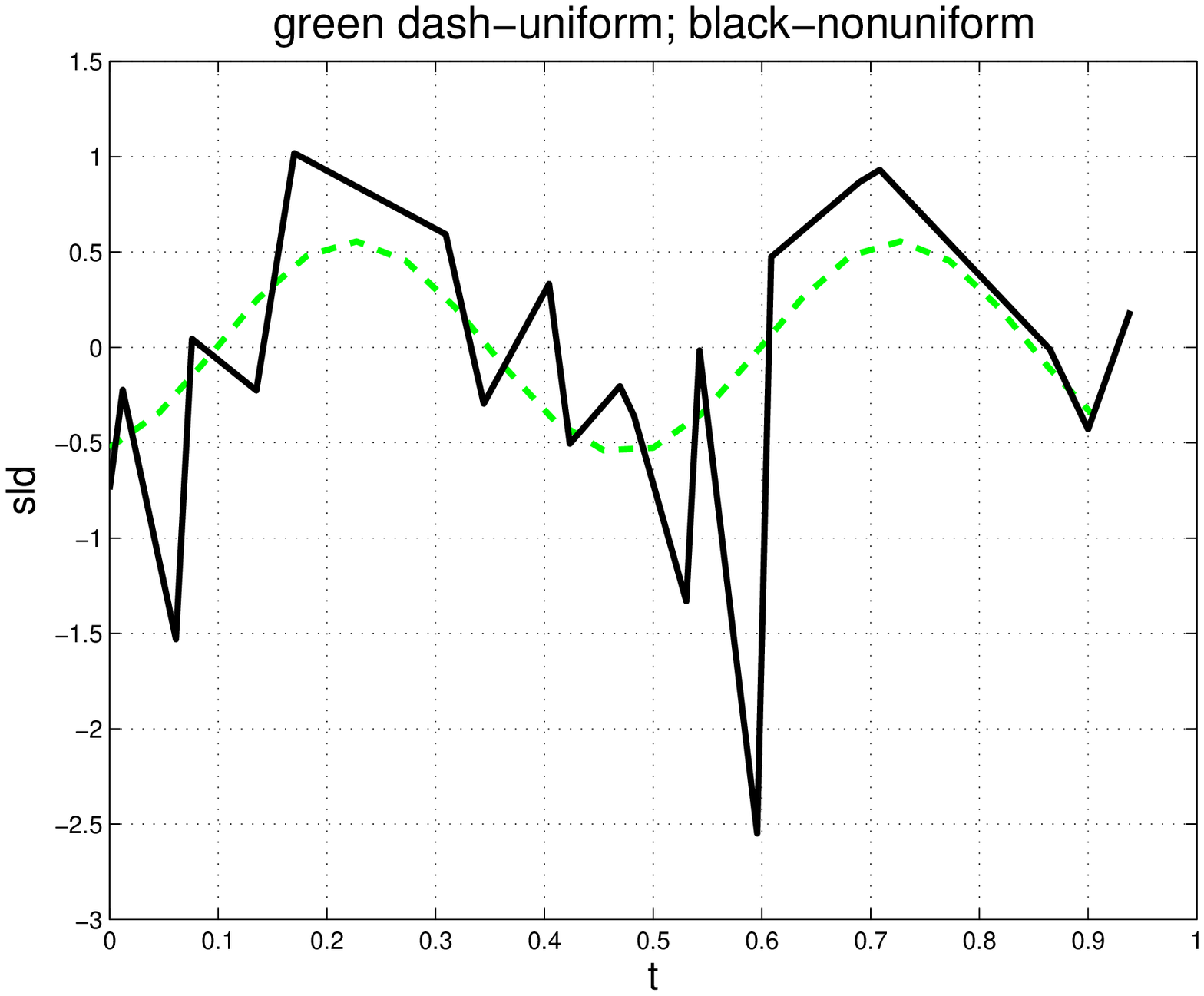,width=1.8in,height=1.6in}

\vspace{3mm}

\parbox[t]{13.5cm}{\footnotesize{\bf Figure 5.3.} Plot of
true solution (red dotted curves) and approximations (blue curves)
on the uniform mesh (LEFT) and the nonuniform mesh (MIDDLE).
Corresponding scaled local difference curves are given (RIGHT: green
for the uniform mesh case and black for the nonuniform case). We may
observe from the first two pictures that, while the approximation on
the uniform mesh is still reasonable, the approximation on the
nonuniform mesh becomes unacceptable. The third picture confirms
these by showing the scaled difference values over the domain (the
black curve acts violently with a significantly large amplitude).}
\end{center}

\vspace{3mm}

\no{\bf Example 5.3.} Consider the second difference $\delta (D_-
f).$ The intention of using the central difference is to improve the
numerical result. It is found in Figure 5.4 that \bbbb
\mbox{sgei}&\approx& 0.2817~~\mbox{ when $\TT$ is uniform}\\
\mbox{sgei}&\approx& 0.7829~~\mbox{ when $\TT$ is nonuniform}\eeee
In addition to a larger global relative error, in Figure 5.4, we may
also observe that if a nonuniform mesh is used, there is a great
irregularity in errors. If the difference is used as an
approximation, very incorrect answers will probably result.

\begin{center}
\epsfig{file=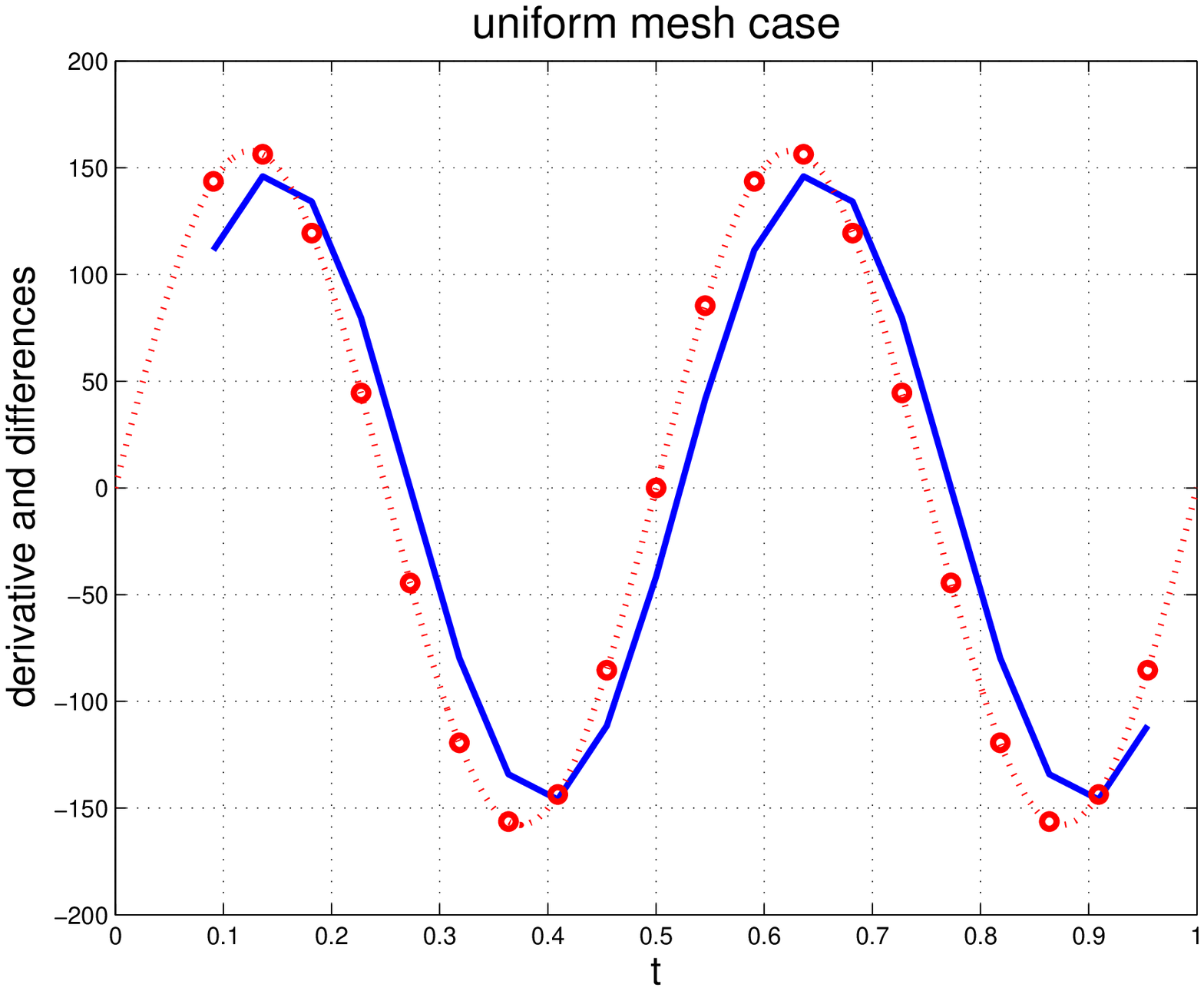,width=1.8in,height=1.6in}
\epsfig{file=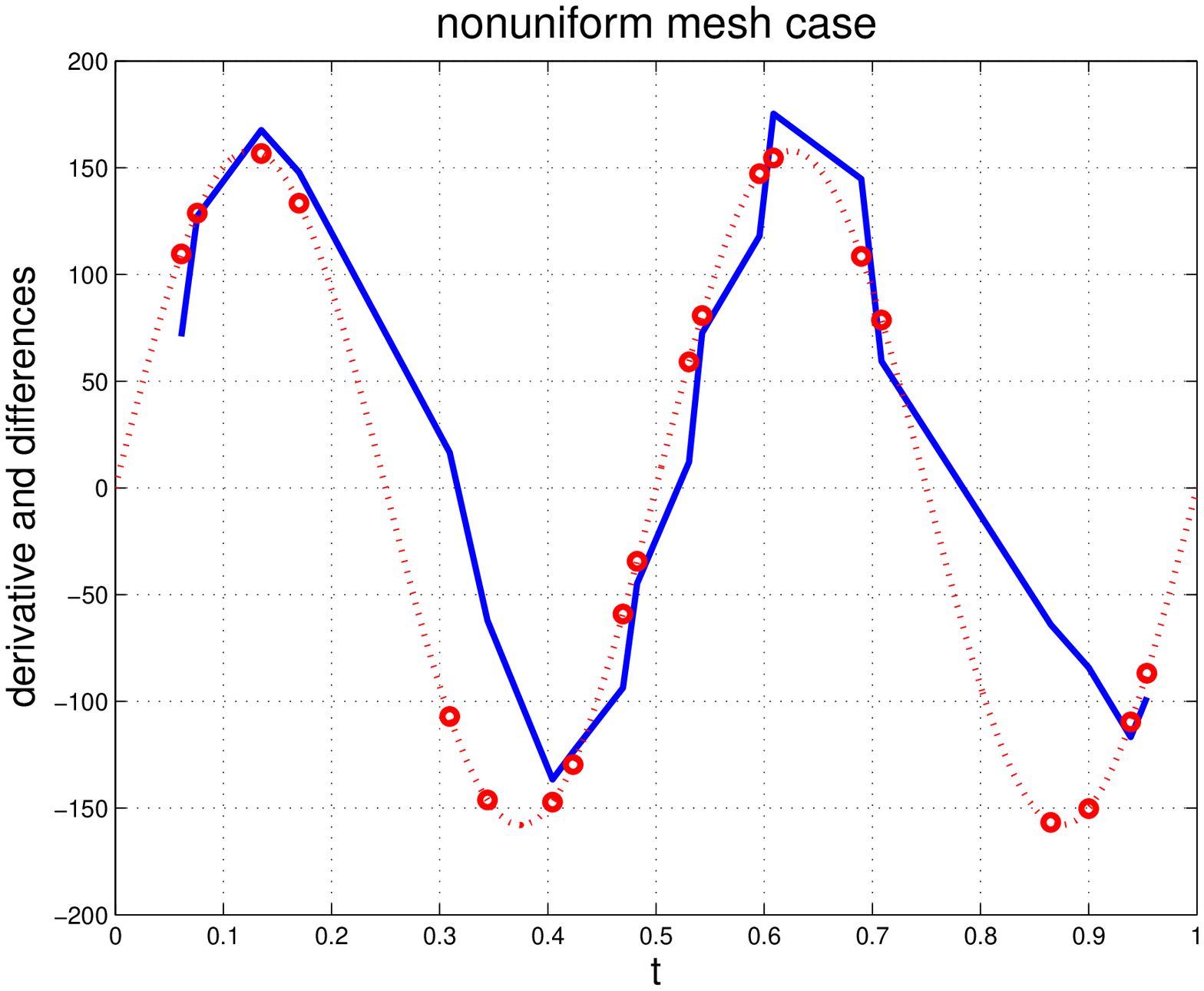,width=1.8in,height=1.6in}
\epsfig{file=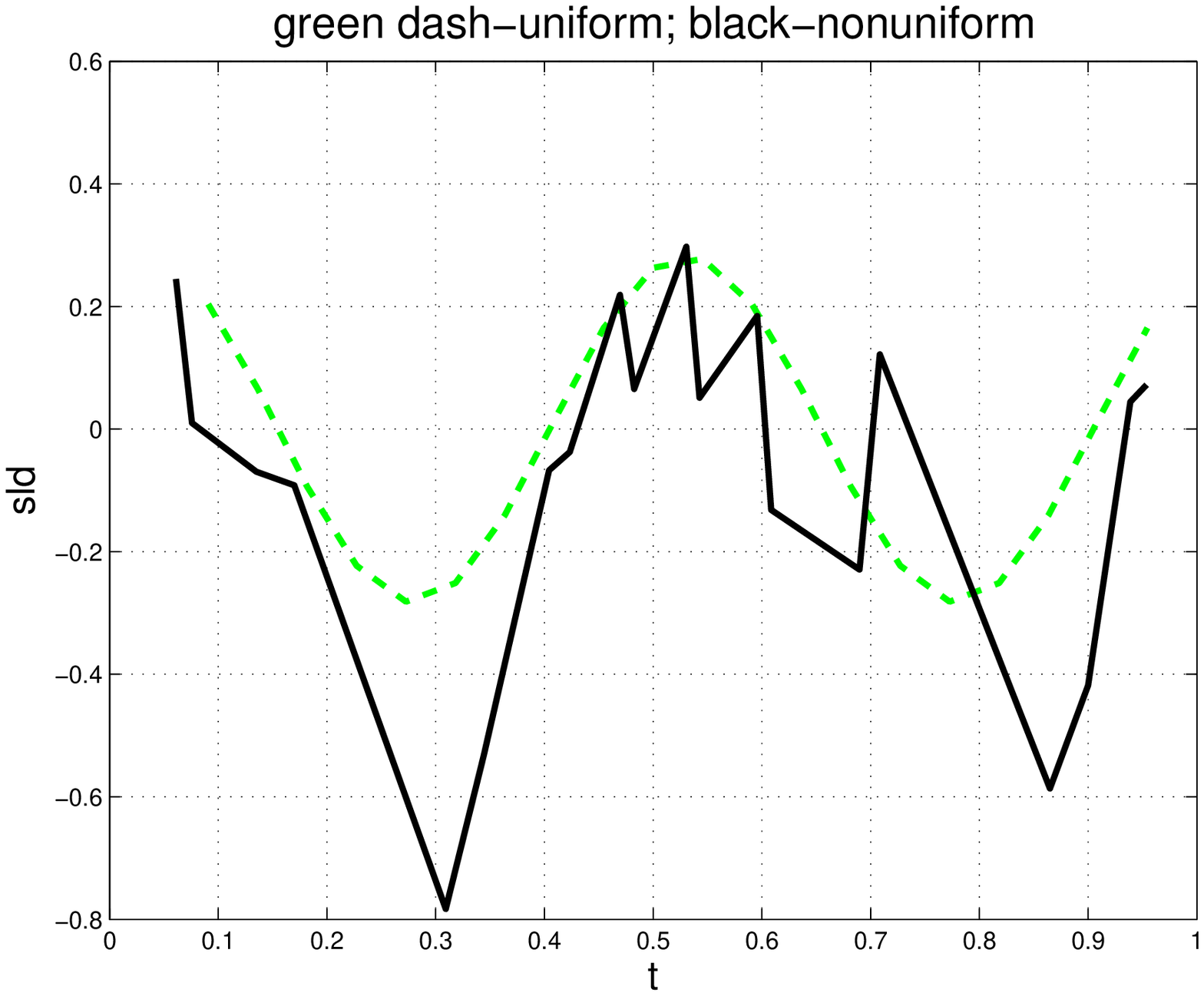,width=1.8in,height=1.6in}

\vspace{3mm}

\parbox[t]{13.5cm}{\footnotesize{\bf Figure 5.4.} Plot of
true solution (red dotted curves) and approximations (blue curves)
on the uniform mesh (LEFT) and the nonuniform mesh (MIDDLE).
Corresponding scaled local difference curves are given (RIGHT: green
for the uniform mesh case and black for the nonuniform case). We may
observe from the first two pictures that, while the approximation on
the uniform mesh looks good, the approximation on the nonuniform
mesh becomes unacceptable. The third picture again confirms these by
showing the scaled difference values over the domain (the black
curve oscillates violently with a large amplitude).}
\end{center}

\vspace{3mm}

\no{\bf Example 5.4.} Consider the second difference $\delta(\delta
f).$ Although the use of central differences continues to improve
the numerical result, it still cannot change the basic features of
the approximations. In this case, we have \bbbb
\mbox{sgei}&\approx& 0.1031~~\mbox{ when $\TT$ is uniform}\\
\mbox{sgei}&\approx& 0.5080~~\mbox{ when $\TT$ is nonuniform}\eeee
The strong irregularity in the nonuniform mesh case shown in Figure
5.5 is not surprising.

\begin{center}
\epsfig{file=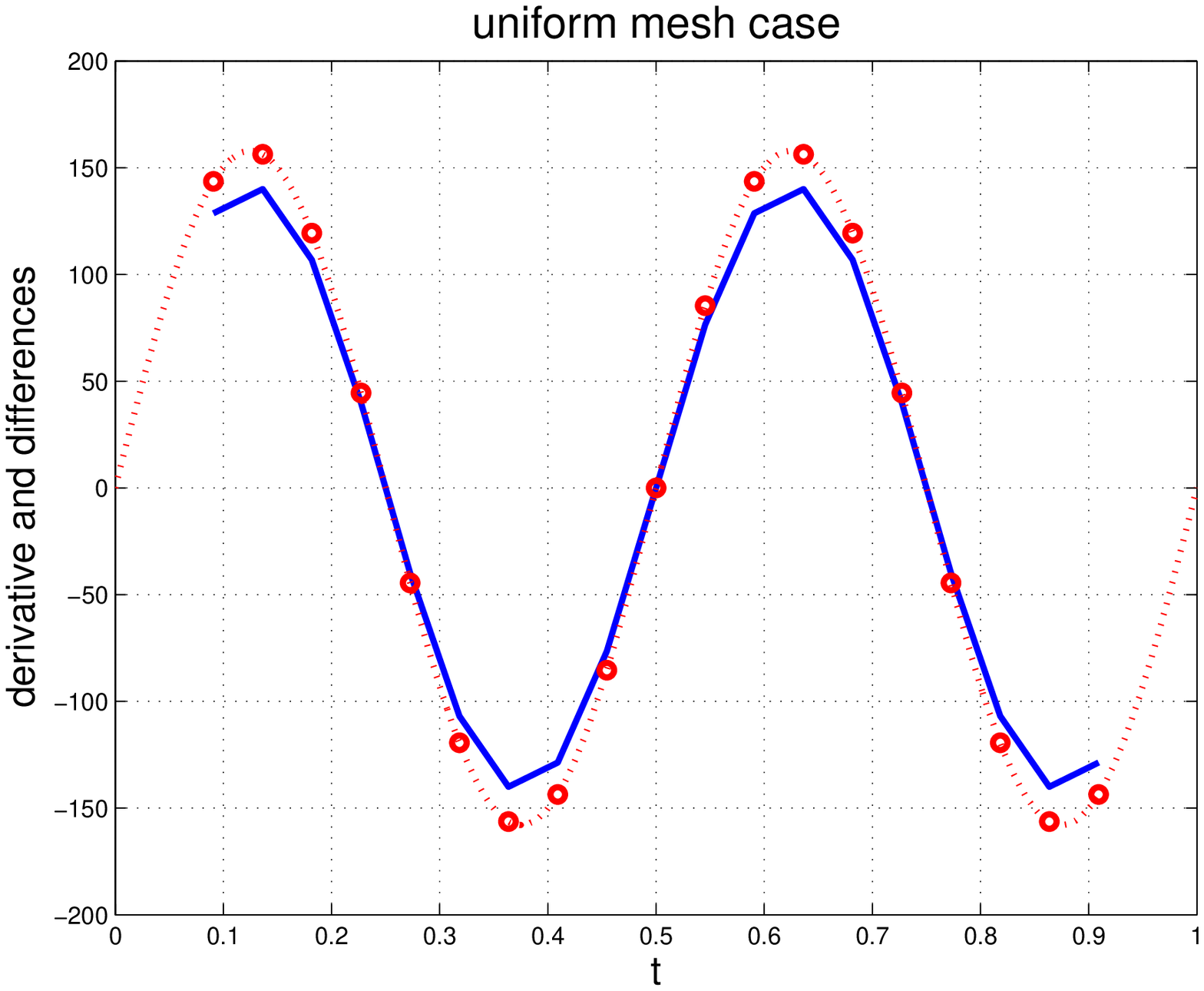,width=1.8in,height=1.6in}
\epsfig{file=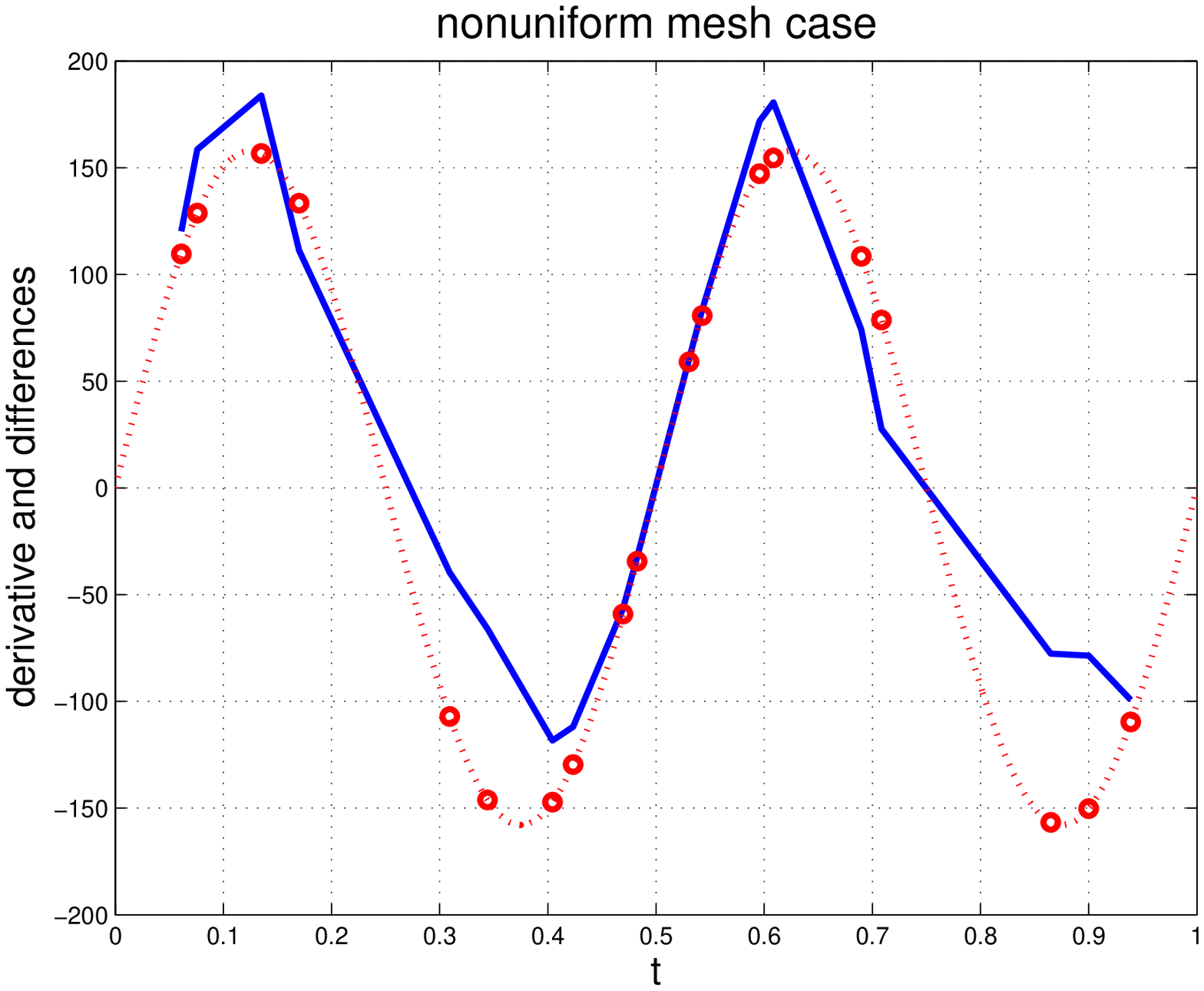,width=1.8in,height=1.6in}
\epsfig{file=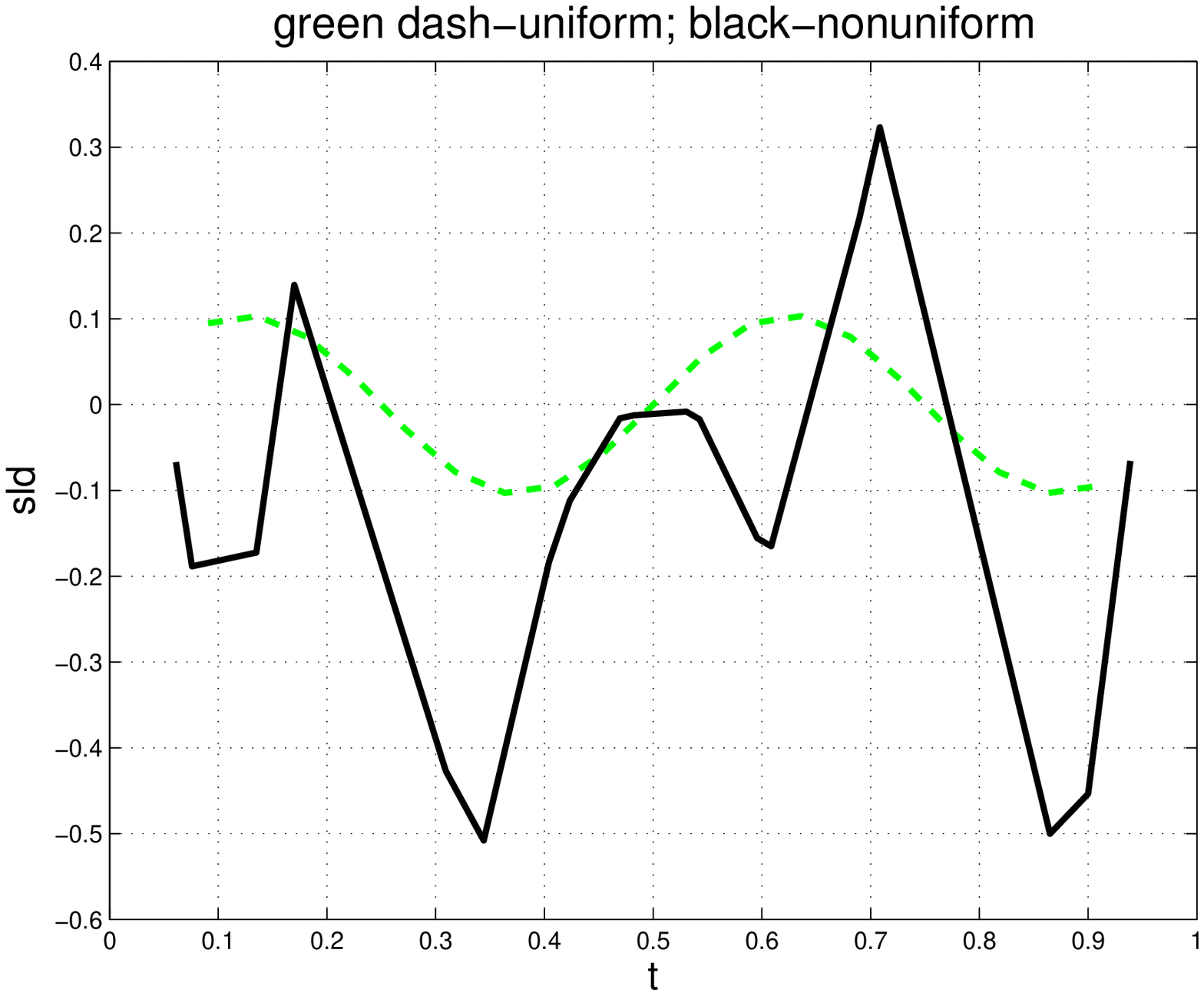,width=1.8in,height=1.6in}

\vspace{3mm}

\parbox[t]{13.5cm}{\footnotesize{\bf Figure 5.5.} Plot of
true solution (red dotted curves) and approximations (blue curves)
on the uniform mesh (LEFT) and the nonuniform mesh (MIDDLE).
Corresponding scaled local difference curves are given (RIGHT: green
for the uniform mesh case and black for the nonuniform case). We may
observe from the first two pictures that, while the approximation on
the uniform mesh is improved due to the use of central difference,
the inconsistency of the formula on the nonuniform mesh remains the
same. These are again confirmed by the third picture through the
scaled difference values over the domain (the black curve oscillates
successively with a large amplitude).}
\end{center}

\vspace{3mm}

\no{\bf Example 5.5.} Let us consider the simple harmonic oscillator
problem where neither a driving force nor friction is assumed
\cite{serway}. If $\phi(t)$ is the displacement of the system at
time $t,$ then the second derivative $\phi''(t)$ is its
acceleration. Based on Hooke's Law and Newton's Second Law, we
obtain the following second order differential equation,
\bb{ee1}\phi''(t)=-\kappa \phi(t),~~~t>t_0,\ee together with the
initial conditions\bb{ee2}\phi(t_0)=1,~~\phi'(t_0)=-1,\ee where
$\kappa$ is a positive constant.

It is not difficult to verify \cite{stewart1} that the solution of
\R{ee1}, \R{ee2} is
\bb{ee3}\phi(t)=-\frac{1}{\sqrt{\kappa}}\sin\sqrt{\kappa}t
+\cos\sqrt{\kappa}t,~~~t\geq t_0.\ee On the other hand, replacing
the second derivative in \R{ee1} by the backward-forward difference,
we acquire \bb{ee4}D_-(D_+w)=-\kappa w(t),~~~t\in\TT,\ee where $\TT$
is a mesh over the interval $[t_0,b].$ We wish the solution of
\R{ee4}, \R{ee2} to be an approximation of the solution of \R{ee1},
\R{ee2}.

First, we let $\TT_1$ be a nonuniform mesh over the interval
$[t_0,b]$ with decreasing steps $h_k=rh_{k-1},~k=1,~2,\ldots,~m,$
for $h_0=1/10,~r=50/59$ and $m=200.$ Second, we let $\TT_2$ be an
uniform mesh with $h=(b-t_0)/10.$

Let $\kappa=4\pi^2,~t_0=0$ and
$b=${\scriptsize$\sum_{k=0}^m$}$h_k\approx 59/90.$ In Figure 5.6
(LEFT), we plot the numerical solutions of \R{ee4}, \R{ee2} on
$\TT_1$ (blue dashed curve), exact solution of the initial value
problem \R{ee1}, \R{ee2} (red dotted curve) together with a
numerical solution of \R{ee4}, \R{ee2} on $\TT_2$ (green curve). A
forward difference is used to approximate the first derivative in
\R{ee2}. Although errors defined for the numerical solution of
initial value problems are slightly different from those used for
function approximations, to illustrate the inconsistency of second
order derivative approximations on nonuniform meshes, let us
continue to use the measurements introduced by Definition 2.4.
Scaled local differences of the numerical solutions on $\TT_1$ (blue
dashed curve) and on $\TT_2$ (green curve) are given in Figure 5.6
(RIGHT).

\vspace{1mm}

\begin{center}
\epsfig{file=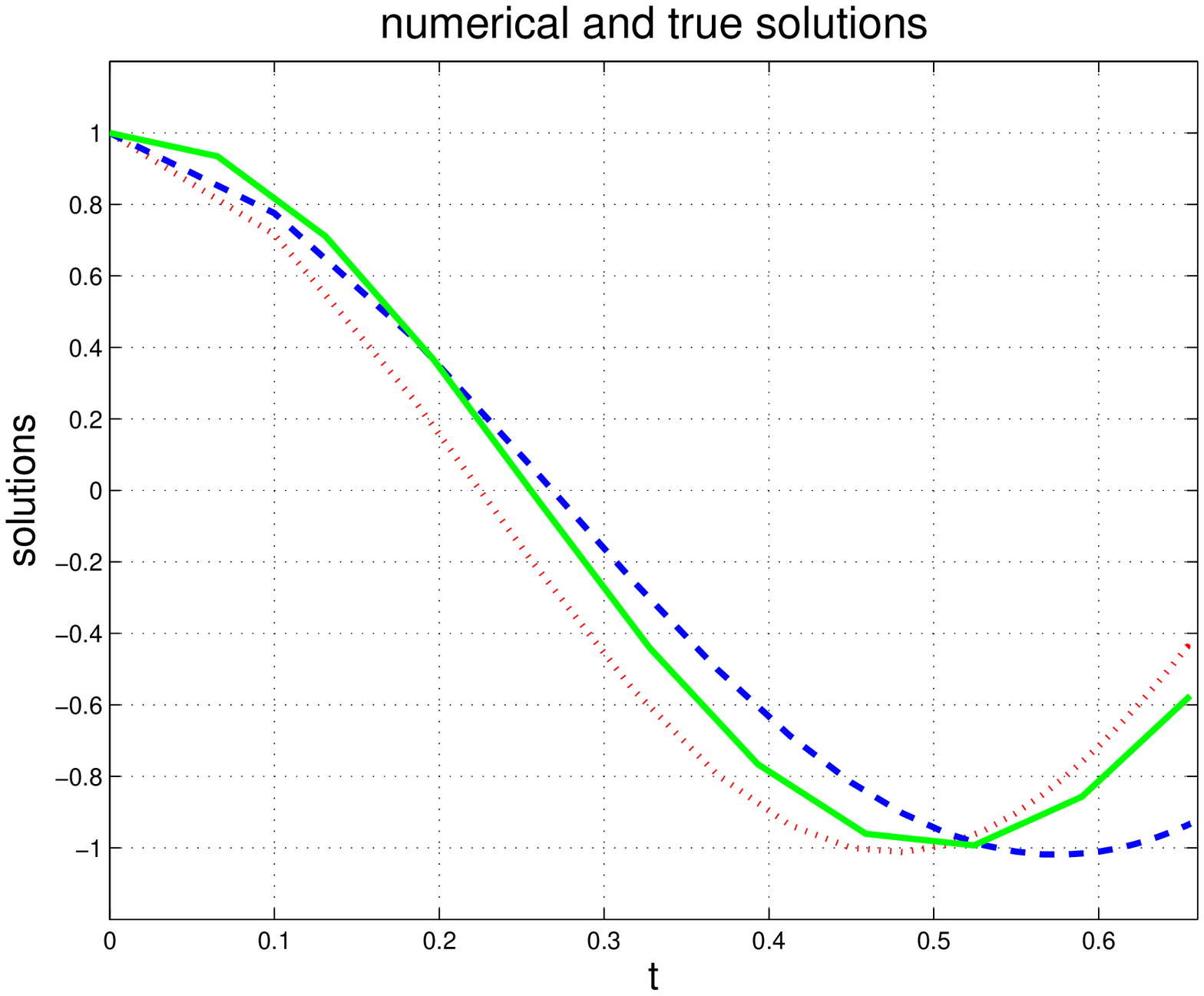,width=2.56in,height=1.6in}
\epsfig{file=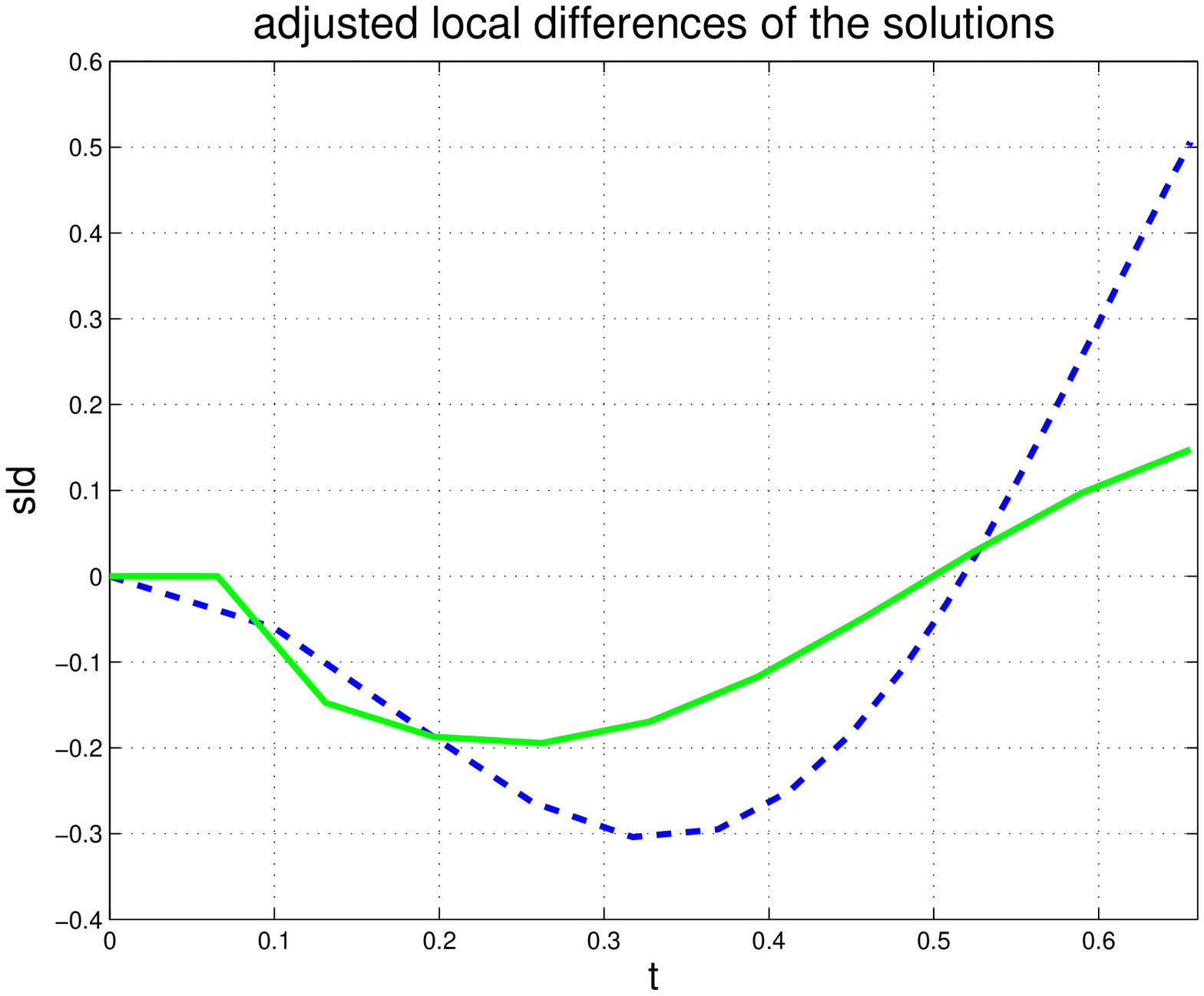,width=2.56in,height=1.6in}

\vspace{2mm}

\parbox[t]{13.5cm}{\footnotesize{\bf Figure 5.6.} Plots of
the solutions of the difference equation problem \R{ee4}, \R{ee2}
and the true solution (LEFT); and the scaled local difference values
of the numerical solutions (RIGHT) on the nonuniform mesh and
uniform mesh, respectively. While the numerical solution and related
sld values on the nonuniform mesh $\TT_1$ are given by blue dashed
curves, the numerical solution and related sld values on the uniform
mesh $\TT_2$ are represented by solid green curves. We may observe
in the first picture that, while the numerical solution on the
uniform mesh maintains a good match to the true solution, which is
indicated by the red dotted curve, the numerical solution on the
nonuniform mesh swings away as $t$ increases. The phenomena are
confirmed by the second picture through the scaled difference values
over the domain used (the blue curve oscillates and its amplitude
increases rapidly as $t\rightarrow b$). }
\end{center}

It is interesting to observe that while the numerical solution
obtained on the uniform mesh keeps a steady error level from the
true solution, although relatively large step is being used
($h=59/900\approx 0.0656$), the numerical solution obtained on the
nonuniform mesh runs away rapidly from the true solution, despite of
the fact that finer steps are employed ($h_k< 0.001$ for $k\geq
14$).

The global error indicator values for the two cases are also
significantly different. It is found that $\mbox{sgei}\approx
0.5055$ on $\TT_1$ and occurs at $x=b$ while $\mbox{sgei}\approx
0.1945$ on $\TT_2$ and occurs at $x\approx 0.2622.$

Since the sizes of $h$ and $h_k$ do not play a major role in the
aforementioned error phenomena, what can be the main cause for the
unsatisfactory approximation on $\TT_1$?

To see the answer, we may notice that
$$\frac{h_k+h_{k-1}}{2h_{k-1}}=\frac{r+1}{2}\neq 1.$$
Therefore, according to \R{yx2} in Theorem 4.2, difference equation
\R{ee4} actually approximates
$$\frac{1+r}{2}\psi''(t)=-\kappa \psi(t),~~~t>t_0,$$ instead of
\R{ee1}! Thus the numerical solution cannot be satisfactory. Of
course, there are several factors, such as stability, that may
affect the numerical solution. Needless to say, however, consistency
is the most fundamental factor.

Conversely, it is not hard to verify that the improved ${\cal
D}_2$ difference given in \R{good} offers a good approximation.
The reader may wish to experiment with the interesting
computations!

\vspace{0.25in}

{\large\bf 6.
Conclusions.~}\setcounter{section}{6}\setcounter{equation}{0} From
the foregoing discussions and numerical experiments, we may
conclude that:
\begin{enumerate}
\item The key issue for any finite difference formula is its
consistency. Only consistent formulas can be used for approximating
derivatives. The effectiveness of an approximation can be measured
by its order of accuracy. For instance, \R{aa2} is an order one
formula according to Theorems 3.1 and 3.2. \item Consistent finite
difference approximations can be derived on uniform or nonuniform
grids. Generally speaking, the higher the order of accuracy is, the
more accurate the finite difference formula can be. \item Difference
approximations on uniform meshes cannot be applied blindly on
nonuniform meshes, nor can difference formulas be composed to form
consistent approximations to second derivatives. At best, they may
lose accuracy; at worst they are not consistent. \item Assuming that
the functions involved are sufficiently smooth, errors of
approximations can be estimated by using Maclaurin series. Errors in
numerical experiments can also be computed via local and global
error formulas. \item The consistency and accuracy of different
approximations can be demonstrated through the use of computer
simulations. Since simulations are based on particular examples,
they are not as rigorous as mathematical proofs.
\end{enumerate}

\vspace{0.20in}

{\large\bf Acknowledgment.~} Many thanks to Dr. Qin Sheng, Professor
of Mathematics at Baylor University, for suggesting the line of
research and for the encouragements and many discussions throughout
this study. The authors are also grateful to Dr. Nancy Miller and
Dr. John Thomason, Professors of Mathematics at Austin Community
College, and Mr. John Miller for reading our revised manuscripts and
for important comments. Last, but not least, the authors would like
to sincerely thank the referee for the many valuable suggestions
which not only helped to improve the content and presentation of
this paper, but also threw light on further study in future
directions.

\end{document}